\documentclass[11pt,a4paper]{article}

\usepackage[T1]{fontenc}
\usepackage{lmodern}
\usepackage[margin=1in]{geometry}
\usepackage{amssymb,mathtools,bm}
\usepackage{amsthm}
\usepackage{graphicx}
\usepackage{enumitem}
\usepackage{microtype}
\usepackage{booktabs}
\usepackage{tikz}
\usepackage{needspace}
\usepackage[hidelinks]{hyperref}
\usepackage[capitalize,nameinlink,noabbrev]{cleveref}

\hypersetup{
  pdftitle={CutFEM for bulk PDE-driven boundary motion described by a level-set function},
  pdfauthor={Erik Burman and Buyang Li}
}
\setlist[itemize]{topsep=3pt,itemsep=2pt,parsep=0pt}
\setlist[enumerate]{topsep=3pt,itemsep=2pt,parsep=0pt,leftmargin=*}

\allowdisplaybreaks[1]
\numberwithin{equation}{section}
\theoremstyle{plain}
\newtheorem{theorem}{Theorem}[section]
\newtheorem{lemma}[theorem]{Lemma}

\newtheorem{proposition}[theorem]{Proposition}
\newtheorem{assumption}[theorem]{Assumption}
\theoremstyle{definition}
\newtheorem{definition}[theorem]{Definition}
\theoremstyle{remark}
\newtheorem{remark}[theorem]{Remark}
\crefname{section}{section}{sections}
\crefname{subsection}{subsection}{subsections}
\Crefname{section}{Section}{Sections}
\Crefname{subsection}{Subsection}{Subsections}
\crefname{assumption}{Assumption}{Assumptions}
\Crefname{assumption}{Assumption}{Assumptions}
\crefname{remark}{Remark}{Remarks}
\Crefname{remark}{Remark}{Remarks}

\newcommand{\R}{\mathbb{R}}
\newcommand{\T}{\mathcal{T}}
\newcommand{\jump}[1]{[\![#1]\!]}
\newcommand{\dd}{\,\mathrm{d}}
\newcommand{\dist}{\operatorname{dist}}

\newcommand{\Om}[1]{\Omega[#1]}
\newcommand{\Ga}[1]{\Gamma[#1]}
\newcommand{\Ub}[1]{U[#1]}
\newcommand{\Oh}{\Omega[\phi_h]}
\newcommand{\Oha}{\Omega_h[\phi_h]}
\newcommand{\Gh}{\Gamma[\phi_h]}
\newcommand{\Uh}{U[\phi_h]}
\newcommand{\Uha}{U_h[\phi_h]}

\newcommand{\Nd}{\mathcal N_{2\delta}}
\newcommand{\nL}[2]{\|#1\|_{L^2(#2)}}
\newcommand{\nH}[2]{\|#1\|_{H^1(#2)}}
\newcommand{\nLi}[2]{\|#1\|_{L^\infty(#2)}}
\newcommand{\nWi}[2]{\|#1\|_{W^{1,\infty}(#2)}}

\title{\bf CutFEM for bulk PDE-driven boundary motion described by a level-set function\thanks{EB was supported in part by EPSRC grants EP/X042650/1 and EP/V050400/1, and BL was supported in part by the NSFC grant 12525111.}}

\author{Erik Burman\thanks{Department of Mathematics, University College
London, Gower Street, London WC1E 6BT, United Kingdom
(\href{mailto:e.burman@ucl.ac.uk}{\texttt{e.burman@ucl.ac.uk}}).}
\quad and\quad Buyang Li\thanks{Department of Applied Mathematics, The Hong Kong
Polytechnic University, Hung Hom, Kowloon, Hong Kong
(\href{mailto:buyang.li@polyu.edu.hk}{\texttt{buyang.li@polyu.edu.hk}}).}}
\date{}


\begin{document}

\maketitle

\begin{abstract}
We analyze a semidiscrete cut finite element method for a quasistatic
reaction--diffusion equation on a domain whose boundary velocity depends on the bulk solution. The method couples a ghost-penalty stabilized bulk solve to transport equations for a level-set variable and its gradient, followed by an $H^1$ reconstruction of the geometry. The transport and reconstruction problems are posed on a band of fixed positive interior width and $O(h)$ exterior width. Ghost penalties on the transport variables and their time derivatives provide control on the active meshes and define the transfer between successive active meshes. The analysis combines level-set shape differentiation, boundary duality, and a transport energy estimate on the exact moving band. A continuation argument closes the coupling between the bulk and geometric errors. For smooth solutions, polynomial degree $k\ge3$, and dimensions two and three, we prove $k$th-order convergence in $H^1$ for the bulk solution and the level set, and in $L^2$ for the transport variables. Numerical experiments with an expanding ellipse show convergence rates consistent with the error estimates, and a five-petal flower that translates, rotates and grows illustrates the effectiveness of the method for rotating domains.
\end{abstract}



\section{Introduction}
\label{sec:intro}

Many free-boundary problems couple a bulk partial differential equation (PDE) to the motion of the domain boundary. Examples
include reaction-controlled dissolution and precipitation \cite{vanNoordenEck2011}, nutrient-dependent tumor growth \cite{FriedmanReitich1999}, and Stefan-type solidification with kinetic undercooling, including quasistatic limits \cite{BarrettGarckeNurnberg2010}. Simulating these processes requires accurate and stable numerical methods for PDE-driven moving boundaries.

A natural approach is to use a moving finite element mesh fitted to the boundary. Arbitrary Lagrangian--Eulerian methods
\cite{BadiaCodina2006,GawlikLew2015}, evolving surface finite element methods (FEMs) \cite{DziukElliott2007}, and the evolving bulk--surface framework of \cite{ElliottRanner2021} provide settings for analyzing numerical approximations on prescribed moving domains and surfaces. Beyond prescribed motion, optimal-order $H^1$ error estimates for evolving FEMs have been established for diffusion-driven surface motion and curvature flows \cite{KLLP2017,KLL2019,KLL2021}, as well as
for a tissue-growth model coupling a bulk Poisson equation to forced mean curvature flow \cite{EdelmannKovacsLubich2025}.
These analyses control the geometric and PDE approximation errors together. However, large deformations may distort a fitted mesh, requiring remeshing and transfer between finite element spaces \cite{GawlikLew2015}. This difficulty motivates unfitted methods, in which the bulk mesh need not conform to the moving boundary.

In unfitted methods the boundary cuts a fixed background mesh. Nitsche-type formulations \cite{HansboHansbo2002}, extended finite element methods \cite{FriesBelytschko2010}, and CutFEM \cite{BurmanEtAl2015} remove the conformity requirement, but must control stability on arbitrarily cut elements. Ghost penalties \cite{Burman2010} and discrete-extension techniques \cite{BurmanHansboLarson2022} address this issue, while high-order formulations \cite{JohanssonLarson2013} and geometric approximations \cite{Lehrenfeld2016,LehrenfeldReusken2018} address the accuracy requirements. The same fixed-background viewpoint underlies surface and narrow-band FEMs \cite{OlshanskiiReuskenGrande2009,DeckelnickEtAl2010,DeckelnickElliottRanner2014}; see \cite{DziukElliott2013} for a survey.

For moving domains, a further difficulty is comparing solutions on successive geometries. Space--time methods treat the evolving integration region directly \cite{HansboLarsonZahedi2016,Zahedi2018}, whereas Eulerian time stepping uses stabilized extension to make previous solutions available on the next domain \cite{LehrenfeldOlshanskii2019}. This approach has been extended to high-order approximation in space and time \cite{LouLehrenfeld2022} and to moving-domain Stokes problems \cite{BurmanFreiMassing2022,vonWahlRichterLehrenfeld2022}. Related analyses treat prescribed evolving surfaces \cite{OlshanskiiReuskenXu2014,LehrenfeldOlshanskiiXu2018}. More recently, Olshanskii and Reusken \cite{OlshanskiiReusken2026} have treated prescribed evolutions with a topological change under structural assumptions near the critical time. In these settings, the geometric approximation is analyzed independently of the computed PDE solution.

This independence is lost when the bulk solution drives the boundary motion. Although unfitted discretizations have been developed for two-phase flow \cite{GrossReusken2007,ClausKerfriden2019,FrachonZahedi2019}, phase change with convection \cite{TchindaNgueyongUrquizaMartin2024a}, solidification \cite{TchindaNgueyongUrquizaMartin2024b}, and fluid--structure interaction \cite{BurmanFernandez2014}, their coupled error analysis must also control the feedback through the computed geometry. Progress in this direction includes the work of von Wahl and Richter \cite{vonWahlRichter2023}, who analyze a parabolic bulk problem coupled to rigid-body dynamics. Their analysis includes the motion error while assuming exact treatment of the geometry once the rigid-body configuration is specified. For a deforming interface, Fern\'andez, Li, and Olshanskii \cite{FernandezLiOlshanskii2026} analyze a semidiscrete CutFEM for fluid--structure interaction, with the interface tracked by a fitted mesh of the elastic solid.

These analyses concern either a rigid body, which has finitely many degrees of freedom, or an elastic solid whose interface is tracked by a fitted mesh. They do not address numerical methods in which the interface is evolved by solving a Hamilton--Jacobi level-set equation coupled to the bulk field. Such level-set formulations are widely used in unfitted computations because they describe the evolving interface implicitly, without requiring a fitted mesh to track it. The additional analytical difficulty is that the level-set discretization error perturbs both the numerical domain and its boundary normal. These perturbations affect the computed bulk solution and, in turn, the velocity driving the subsequent interface evolution.

In this paper we construct and analyze a CutFEM in this setting. We consider the quasistatic problem
\begin{align}
 u-\Delta u=f(\cdot,t)\ \text{in }\Omega(t),\qquad
 \partial_nu=0\ \text{on }\Gamma(t)=\partial\Omega(t),
 \label{eq:intro-model}
\end{align}
where $\Gamma(t)$ moves with outward normal velocity $V_n=\beta(u)$ and $\beta\colon\R\to\R$ is a given smooth function. We represent the evolving domain by a level-set function $\phi$ satisfying
\begin{align}
 \partial_t\phi+\beta(u)|\nabla\phi|=0
 \label{eq:intro-model-phi}
\end{align}
in a band region containing the boundary $\Gamma(t)$; see \cref{subsec:meshes}. This model provides a simple setting in which to study the two-way coupling between an elliptic bulk problem and a numerically evolved level-set geometry. Related PDE-driven free-boundary models arise in tissue growth; see, for example, \cite{EylesKingStyles2019}, where the interface evolves by forced mean curvature flow with forcing determined by a PDE in the enclosed bulk domain. 

The construction of the method is guided by the estimates needed to close this feedback. We evolve the auxiliary variables
$p=\phi$ and $\eta=\nabla\phi$ through transport equations and reconstruct the geometry from $\phi-\Delta\phi=p-\nabla\cdot\eta$, with the weak boundary realization specified in \cref{subsec:scheme}. The reconstruction of $\phi$ controls the level-set $H^1$ error by the transport $L^2$ errors, and the level-set $H^1$ error controls the bulk $H^1$ error without a negative power of $h$ in the coupling; see \cref{prop:eu-final}. The one-sided band has $O(1)$ interior width to provide the required uniform trace bound, and exterior width $O(h)$ so that the bulk ghost extension supplies the transport coefficients. 

To estimate the transport errors, we measure their energy on the exact moving band, where the boundary terms from Reynolds'
formula cancel the transport fluxes. Replacing the numerical band by the exact band introduces a geometric defect involving
the transport residuals. Ghost penalties on the time derivatives control these residuals on the active mesh and make the
defect estimable in the same error norms; see \cref{lem:dt-transport,lem:geometry-defect}.

The finite element functions live on active meshes that are rebuilt only when needed. At an update time $t_m$ we keep the elements near the band,
$$
 U_h=\bigcup\bigl\{K\in\T_h:\ \dist\bigl(K,U[\phi_h(t_m^-)]\bigr)\le c_{\rm buf}h\bigr\},
$$
and keep $U_h$ fixed as long as the moving band stays inside it with a margin,
$$
 \bigl\{\dist(\cdot,U[\phi_h(t)])\le\tfrac12c_{\rm buf}h\bigr\}\Subset U_h
 \subset\bigl\{\dist(\cdot,U[\phi_h(t)])\le 2c_{\rm buf}h\bigr\}.
$$
When this fails at $t_{m+1}$, a new $U_h$ is built and the transport variables are moved to it by a stabilized extension,
$$
 p_h(t_{m+1}^+)=E_{U_h}\,p_h(t_{m+1}^-),\qquad
 \eta_h(t_{m+1}^+)=E_{U_h}\,\eta_h(t_{m+1}^-);
$$
see \cref{subsec:active,subsec:scheme}. The active mesh thus moves with the boundary without remeshing.

Finally, the transfer at each update has an error amplification factor close to one; see \cref{lem:transfer}. A lower bound $ch$ for the time interval between active-mesh updates bounds the number of updates by $C/h$. Combining the transfer and fixed-mesh estimates gives $k$th-order convergence in the $H^1$ norm for $u$ and $\phi$ (and in the $L^2$ norm for $p$ and $\eta$) on their respective active meshes, for smooth evolution in dimensions two and three and finite elements of polynomial degree $k\ge3$; see \cref{thm:main}. 

The rest of this paper is organized as follows. \Cref{sec:continuous,sec:method,sec:main} present the problem, the method and the main result. \Cref{sec:prelim,sec:error-eq,sec:apriori} contain the proof, \cref{sec:num} the numerical experiments, and \cref{sec:conclusion} the conclusion. The appendices give the full discretization and the remaining parts of the proof.

\section{The continuous problem}
\label{sec:continuous}

Let $D\subset\R^n$, $n\in\{2,3\}$, be a bounded Lipschitz hold-all domain, and let $\Omega(t)\Subset D$, $t\in[0,T]$, have boundary $\Gamma(t)=\partial\Omega(t)$. The boundary moves with outward normal velocity
\begin{equation}
 V_n=\beta(u|_{\Gamma(t)}),
 \label{eq:normal-speed}
\end{equation}
where $\beta\colon\R\to\R$ is a given function, and $u$ is a bulk field in $\Omega(t)$ determined by the elliptic problem
\begin{equation}
 u-\Delta u=f(\cdot,t)\ \text{ in }\Omega(t),\qquad
 \partial_nu=0\ \text{ on }\Gamma(t) .
 \label{eq:bulk}
\end{equation}
Here $f\colon D\times[0,T]\to\R$ is a given source. 

\subsection{Level-set description}
We describe the moving boundary by a level-set function $\phi\colon D\times[0,T]\to\R$: 
\begin{subequations}\label{eq:levelset-domain}
\begin{equation}
 \Gamma(t)=\{x\in D:\phi(x,t)=0\} ,
\end{equation}
with 
\begin{align}
&\{x\in D:-2\delta<\phi(x,t)<0\}\subset\Omega(t) ,\\
&\{x\in D:0<\phi(x,t)<2\delta\}\cap\overline{\Omega(t)} = \emptyset,
\end{align}
\end{subequations}
where $\delta>0$ is a prescribed constant independent of $t\in[0,T]$. Moreover, $\phi$ is nondegenerate near the boundary: 
\begin{equation}
 |\nabla\phi|\ge c>0 \quad\mbox{on}\,\,\, \Nd(t):=\{x\in D:|\phi(x,t)|<2\delta\} .
 \label{eq:Nd-def}
\end{equation}
We choose the sign of $\phi$ so that $\Omega(t)=\{x\in D:\phi(x,t)<0\}$,
and assume that $\inf_{t\in[0,T]}\dist(\overline{\Nd(t)},\partial D)>0$.

Since the velocity is needed on the whole band $\Nd(t)$, which extends
outside $\Omega(t)$, we fix a smooth extension of $u$ from $\Omega(t)$ to
$D$, still denoted by $u$, and let $\phi$ evolve by
\begin{equation}
 \partial_t\phi+\beta(u)n(\nabla\phi)\cdot\nabla\phi=0
 \quad\text{in } \Nd(t),\ t\in(0,T],
 \qquad n(\nabla\phi) = \frac{\nabla\phi}{|\nabla\phi|}.
 \label{eq:HJ}
\end{equation}
The initial value is $\phi(\cdot,0)=\phi_0$. The extension affects $\phi$ away from $\Gamma(t)$ but not the motion of
$\Gamma(t)$, whose normal velocity is $\beta(u|_{\Gamma(t)})$.

\subsection{Reformulation of the problem}
We use the augmented variables
\begin{equation}
 p:=\phi,\qquad
 \eta:=\nabla\phi=(\eta_1,\ldots,\eta_n)^{\rm T},\qquad
 \eta_j:=\partial_j\phi ,
 \label{eq:peta-def}
\end{equation}
which satisfy the following evolution equations:
\begin{align}
 \partial_tp+\beta(u)n(\nabla\phi)\cdot\nabla p&=0
 &&\text{in }\Nd(t),
 \label{eq:p-cont}\\
 \partial_t\eta_j+\beta(u)n(\nabla\phi)\cdot\nabla\eta_j
 &=-\beta'(u)\partial_ju\,|\nabla\phi|
 &&\text{in }\Nd(t),\quad j=1,\ldots,n .
 \label{eq:eta-cont}
\end{align}
Equation \eqref{eq:p-cont} is \eqref{eq:HJ}; differentiation in $x_j$ gives \eqref{eq:eta-cont}, since
$$
 \partial_j|\nabla\phi|
 =n(\nabla\phi)\cdot\partial_j\nabla\phi
 =n(\nabla\phi)\cdot\nabla\eta_j.
$$
We evolve $p,\eta$ by \eqref{eq:p-cont}--\eqref{eq:eta-cont} and reconstruct $\phi$ from
\begin{align}
\phi-\Delta\phi = p - \nabla\cdot\eta 
 \label{eq:phi-reconstr}
\end{align}
on a numerically computed band, with the natural Neumann condition $(\nabla\phi-\eta)\cdot n(\nabla\phi)=0$ on the boundary of the band. 

\section{The numerical method}
\label{sec:method}

We define the bands, the active meshes and the ghost penalties, and then the semidiscrete CutFEM. A fully discrete scheme is given in \cref{subsec:fd-scheme}.

\subsection{Background mesh and level-set bands}
\label{subsec:meshes}

Let $\T_h$ be a quasiuniform, shape-regular simplicial mesh of $D$ with mesh size $h$, and let
\begin{equation}
 V_h:=\{v_h\in C^0(\bar D):v_h|_K\in\mathbb P_k(K)\ \forall K\in\T_h\},
 \qquad k\ge1 ,
 \label{eq:Vh}
\end{equation}
be the space of continuous piecewise polynomials of degree $k$. For $\chi=\phi$ or the discrete level set $\chi=\phi_h$, extended to $D$ as in \cref{subsec:active}, and $c\in\R$, we write
\begin{align}\label{eq:Omega-Gamma-psi}
\begin{aligned}
 \Ga{\chi}&:=\{\chi=0\} 
 \quad\mbox{and}\quad
 \Gamma_c[\chi]:=\{\chi=c\}, \\
 \Om{\chi}&:=\{x\in D:\chi(x)<0\}.
\end{aligned}
\end{align} 

Let $\kappa\geq 1$ be a fixed constant independent of $h$.  The \emph{one-sided band} of $\chi$ and its two boundaries are defined as 
\begin{equation}
 \Ub{\chi}:=\{x\in D:-\delta<\chi(x)<\kappa h\},\quad
 \Gamma_+[\chi]:=\Gamma_{\kappa h}[\chi],\quad
 \Gamma_-[\chi]:=\Gamma_{-\delta}[\chi].
 \label{eq:Gamma-c-def}
\end{equation}
When $\chi$ is nondegenerate near these level surfaces, $\Ub{\chi}$ is a neighborhood of $\Ga{\chi}$ of thickness $\simeq h$ on the exterior side and $\simeq\delta$ on the interior side. The discrete level-set function $\phi_h$ (defined below) determines the
numerical bulk domain $\Oh$, the numerical interface $\Gh$, and the numerical band $\Uh$; the exact band is $\Ub{\phi}$.

\subsection{Active meshes}
\label{subsec:active}

The finite element functions are defined on active meshes. These are rebuilt at times $0=t_0<t_1<\dots<t_M=T$ and kept fixed on each interval $I_m=(t_m,t_{m+1}]$. For a function $w$ of time, $w(t_m^-)$ denotes its left limit at $t_m$, computed on the old meshes, and $w(t_m^+)$ its value after the update at $t_m$, on the new meshes. At $t_m$ the active meshes are built from $\phi_h(t_m^-)$, with $\phi_h(t_0^-):=I_h\phi_0$. For a fixed constant $c_{\rm buf}\ge4$ (as a convenient sufficient choice), let
\begin{align}
 \mathcal B_{h,m}
 &:=
 \{K\in\T_h:
   \dist(K,\Ub{\phi_h(t_m^-)})\le c_{\rm buf}h\},
 \label{eq:buffer-elements}\\
 \T_{h,m}^\Omega
 &:=
 \{K\in\T_h:
   K\cap\Omega[\phi_h(t_m^-)]\ne\varnothing\}
 \cup\mathcal B_{h,m},\\
 \Oha
 &:=
 \bigcup_{K\in\T_{h,m}^\Omega}K,
 \qquad
 \Uha:=
 \bigcup_{K\in\mathcal B_{h,m}}K
 \subset\Oha
 \quad\text{for }t\in I_m.
 \label{eq:numerical-bulk-support}
\end{align}
Thus $\Uha$ consists of the elements within distance $c_{\rm buf}h$ of the band, and $\Oha$ adds the elements that meet the bulk domain; see \cref{fig:active-supports-schematic}. We write $\Oha$ and $\Uha$ without explicitly indicating the update index $m$, but keep in mind that these active meshes are constructed from the pre-update level set $\phi_h(t_m^-)$ and remain fixed on $I_m$; they are not determined by the current $\phi_h(t)$.

The level set $\phi_h$ is computed only on $\Uha$. To define $\Omega[\phi_h]$ and $\Ub{\phi_h}$ on all of $D$, we set $\phi_h:=-2\delta$ on the connected components of $D\setminus\Uha$ that lie in $\Omega[\phi_h(t_m^-)]$, and $\phi_h:=\kappa h+\delta$ on the other components. Since $\Ga{\phi_h(t_m^-)}\subset\Uha$, each component lies either in $\Omega[\phi_h(t_m^-)]$ or outside it. These values are not in $[-\delta,\kappa h]$, so $\Ub{\phi_h(t)}\subset\Uha$. This extension is used only to describe the geometry.

\paragraph{Update times}
We keep the meshes as long as
\begin{subequations}\label{eq:band-buffer}
\begin{align}
 \bigl\{x\in D:
        \dist(x,\Uh)\le\tfrac12 c_{\rm buf}h
 \bigr\}
 &\Subset\Uha, \label{eq:band-buffer-a}\\
 \bigl\{x\in D:
        \dist(x,\Uh)\le 2c_{\rm buf}h
 \bigr\}
 &\supset\Uha, \label{eq:band-buffer-b}
\end{align}
\end{subequations}
that is, the band stays more than $\tfrac12c_{\rm buf}h$ away from $\partial\Uha$, and $\Uha$ stays within $2c_{\rm buf}h$ of the band. Let
$$
 t_{m+1}:=\min\bigl\{T,\ \inf\{t>t_m:\ \eqref{eq:band-buffer}\text{ fails}\}\bigr\}.
$$
Then \eqref{eq:band-buffer} holds for $t_m<t<t_{m+1}$, and its non-strict version holds for $\phi_h(t_{m+1}^-)$. If $t_{m+1}<T$, new meshes are built from $\phi_h(t_{m+1}^-)$ and the solution is transferred to them as in \cref{subsec:scheme}. Under the assumptions of \cref{thm:main} and for $h$ small enough, the update times satisfy
\begin{equation}
 t_{m+1}-t_m\ge c_{\rm interval}h
 \quad\text{whenever }t_{m+1}<T ,
 \label{eq:event-time interval-spacing}
\end{equation}
with $c_{\rm interval}>0$ independent of $h$; see \cref{subsec:interval}. Hence there are at most $C/h$ updates.

\input{figures_m2/active_supports_figure.tex}

\subsection{Ghost penalties and the stabilized extension}
\label{subsec:ghost}

Let $\mathcal F_h^g(\Uha)$ be the set of interior faces of $\Uha$, and consider the face-based ghost penalties of \cite{Burman2010} with two scalings $r\in\{0,1\}$:
\begin{equation}
 g_r(w_h,v_h)_{\Uha}
 :=\sum_{F\in\mathcal F_h^g(\Uha)}\sum_{\ell=1}^k
 \gamma_{r,\ell}\, h^{2\ell+1-2r}
 \bigl(\jump{\partial_{n_F}^\ell w_h},
       \jump{\partial_{n_F}^\ell v_h}\bigr)_F ,
\label{eq:ghost-bulk-general}
\end{equation}
where $\jump{\cdot}$ denotes the jump across the face $F$, $n_F$ is one fixed unit normal used for both traces, $\partial_{n_F}^\ell$ is the $\ell$th derivative in that direction, and $\gamma_{r,\ell}>0$ are parameters independent of $h$; we write $\gamma_0:=\min_\ell\gamma_{0,\ell}$. The ghost penalties give the extension estimate of \cref{lem:ghost-extension}, whose constants do not depend on the cutting position.

\begin{definition}[Chain condition]
\label{def:chain}
Let $\mathcal S_h\subset\T_h$ and $Q_h\subset D$. The pair $(\mathcal S_h,Q_h)$ satisfies the chain condition with constant $N_*$
if every $K\in\mathcal S_h$ with $K\not\subset Q_h$ admits a chain of distinct elements $\mathcal P_K=(K_0,\ldots,K_\ell)$ such that
$$
 \begin{aligned}
  K_0&=K,\qquad K_\ell\subset Q_h,\qquad
  K_i\in\mathcal S_h
  &&(0\le i\le\ell),\\
  K_{i-1}\cap K_i&\in\mathcal F_h^g(\Uha)
  &&(1\le i\le\ell),
 \end{aligned}
$$
and the chains have bounded length and overlap:
$$
 \ell+1\le N_*,
 \qquad
 \#\{\mathcal P_K:L\in\mathcal P_K\}\le N_*
 \quad\text{for every }L\in\mathcal S_h.
$$
\end{definition}

The chain condition is verified in \cref{section:patch,appendix:path} for $t<t^\ast$ and $h$ small enough.

\begin{lemma}[Cut-independent ghost extension {\cite[Lemma~5.5]{LehrenfeldOlshanskii2019}}]
\label{lem:ghost-extension}
Suppose that $(\mathcal B_{h,m},\Uh)$ and $(\T_{h,m}^\Omega,\Oh)$ satisfy the chain condition with a constant $N_*$ independent of $h$, the cut position, and the update index. Then, for every $v_h\in V_h$,
\begin{equation}
\begin{aligned}
 \nL{v_h}{\Uha}^2
 &\le C_{\rm ext}\bigl(
       \nL{v_h}{\Uh}^2
       +\gamma_0^{-1}g_0(v_h,v_h)_{\Uha}\bigr),\\
 \nH{v_h}{\Uha}^2
 &\le C_{\rm ext}\bigl(
       \nH{v_h}{\Uh}^2
       +g_1(v_h,v_h)_{\Uha}\bigr),\\
 \nH{v_h}{\Oha}^2
 &\le C_{\rm ext}\bigl(
       \nH{v_h}{\Oh}^2
       +g_1(v_h,v_h)_{\Uha}\bigr).
\end{aligned}
\label{eq:proved-ghost-extension}
\end{equation}

The constant $C_{\rm ext}$ is independent of $h$, the cut position, and the update index.
\end{lemma}

Let $U_h$ be an active mesh and $U\subset U_h$ a band such that $(\mathcal B_{h,m},U)$ satisfies the chain condition. The stabilized extension $E_{U_h}w\in V_h|_{U_h}$ of $w\in L^2(U)$ is defined by
\begin{equation}
 (E_{U_h}w,v_h)_{U}+g_0(E_{U_h}w,v_h)_{U_h}=(w,v_h)_{U}
 \qquad\forall v_h\in V_h|_{U_h} .
 \label{eq:stabilized-extension}
\end{equation}
By the first estimate in \eqref{eq:proved-ghost-extension} with $U$ in place of $\Uh$, the left-hand side is coercive on $V_h|_{U_h}$, so $E_{U_h}w$ is well defined. Taking $v_h=E_{U_h}w$ gives
\begin{equation}
 \nL{E_{U_h}w}{U}^2+g_0(E_{U_h}w,E_{U_h}w)_{U_h}
 \le\nL{w}{U}^2.
 \label{eq:stabilized-extension-stability}
\end{equation}

\subsection{Semidiscrete CutFEM with active-mesh updates}
\label{subsec:scheme}

For a union $S$ of elements, let $V_h|_S:=\{v_h|_S:v_h\in V_h\}$. On each interval $I_m$ we seek $u_h\in V_h|_{\Oha}$, $p_h,\phi_h\in V_h|_{\Uha}$ and $\eta_h=(\eta_{1,h},\ldots,\eta_{n,h})^{\rm T}\in (V_h|_{\Uha})^n$ such that (i)--(iii) below hold.

\paragraph{(i) Transport of $p_h$ and $\eta_h$} For all $v_h\in V_h|_{\Uha}$ and $j=1,\dots,n$,
\begin{align}
 &(\partial_t p_h+\beta(u_h)n(\eta_h)\cdot\nabla p_h,v_h)_{\Uh}
 +g_0(\partial_tp_h+p_h,v_h)_{\Uha}=0,                                \label{eq:scheme-p}\\
 &(\partial_t \eta_{j,h}+\beta(u_h)n(\eta_h)\cdot\nabla\eta_{j,h},v_h)_{\Uh}
 +g_0(\partial_t\eta_{j,h}+\eta_{j,h},v_h)_{\Uha}\nonumber\\
 &\qquad=-(\beta'(u_h)\,\partial_ju_h\,|\eta_h|,v_h)_{\Uh} ,    \label{eq:scheme-eta}
\end{align}
where $n(\eta_h)=\eta_h/|\eta_h|$. We use $n(\eta_h)$ rather than $n(\nabla\phi_h)$ because $\nabla\phi_h$ is discontinuous across element faces. For $t<t^\ast$, \eqref{eq:smallness}, \eqref{eq:Nd-def} and $\Uh\subset\Nd(t)$ (see \cref{subsec:continuation}) imply that $|\eta_h|$ is bounded away from zero on $\Uh$, so $n(\eta_h)$ is well defined.

\paragraph{(ii) Level-set reconstruction}
Given $(p_h,\eta_h)$, find $\phi_h\in V_h|_{\Uha}$ such that
\begin{equation}
\begin{aligned}
 &(\phi_h,v_h)_{\Uh}+(\nabla\phi_h,\nabla v_h)_{\Uh}
 +g_1(\phi_h,v_h)_{\Uha}\\
 &\qquad=(p_h,v_h)_{\Uh}+(\eta_h,\nabla v_h)_{\Uh}
\end{aligned}
\label{eq:scheme-rec}
\end{equation}
for all $v_h\in V_h|_{\Uha}$. Since $\Uh$ depends on $\phi_h$, this problem is nonlinear. When $(p_h,\eta_h)$ is close to $(p,\eta)$, it has a unique solution near the exact level set, which we call the local branch; see \cref{sec:ODE}.

\paragraph{(iii) Bulk solve}
Given $\phi_h$, find $u_h\in V_h|_{\Oha}$ such that
\begin{equation}
 (u_h,v_h)_{\Oh}+(\nabla u_h,\nabla v_h)_{\Oh}
 +g_1(u_h,v_h)_{\Uha}
 =(f(t),v_h)_{\Oh}
\label{eq:scheme-u}
\end{equation}
for all $v_h\in V_h|_{\Oha}$.

Eliminating $\phi_h$ and $u_h$ through (ii) and (iii) turns (i) into a system of ordinary differential equations for the nodal values of $(p_h,\eta_h)$ on each $I_m$; see \cref{sec:ODE}.

\paragraph{Initial data}
With $p(0)=\phi_0$ and $\eta(0)=\nabla\phi_0$, we set
\begin{equation}
 p_h(0)=I_hp(0),\qquad \eta_h(0)=I_h\eta(0),
 \label{eq:interpolated-initial-data}
\end{equation}
where $I_h$ is the nodal interpolation operator onto $V_h$. The initial active meshes are built from $I_h\phi_0$ as in \cref{subsec:active}. Then $\phi_h(0)$ and $u_h(0)$ are given by \eqref{eq:scheme-rec} and \eqref{eq:scheme-u}.

\paragraph{Transfer at update times}
At an update time $t_{m+1}<T$ (see \cref{subsec:active}), the transport variables are transferred to the new active mesh $U_h[\phi_h(t_{m+1}^-)]$ by the stabilized extension \eqref{eq:stabilized-extension}:
\begin{equation}
 p_h(t_{m+1}^+):=E_{U_h[\phi_h(t_{m+1}^-)]}\,p_h(t_{m+1}^-),\qquad
 \eta_h(t_{m+1}^+):=E_{U_h[\phi_h(t_{m+1}^-)]}\,\eta_h(t_{m+1}^-).
 \label{eq:transfer}
\end{equation}
Both extensions use the band $U[\phi_h(t_{m+1}^-)]$. By the non-strict version of \eqref{eq:band-buffer-a}, this band lies in the old mesh $U_h[\phi_h(t_m^-)]$, so $p_h(t_{m+1}^-)$ and $\eta_h(t_{m+1}^-)$ are defined on it; by \eqref{eq:buffer-elements}, it also lies in the new mesh $U_h[\phi_h(t_{m+1}^-)]$. 

Next, $\phi_h(t_{m+1}^+)\in V_h|_{U_h[\phi_h(t_{m+1}^-)]}$ is computed from \eqref{eq:scheme-rec} with the transferred data
$p_h(t_{m+1}^+),\eta_h(t_{m+1}^+)$. In \eqref{eq:scheme-rec} the integration band is $U[\phi_h(t_{m+1}^+)]$. Then 
$u_h(t_{m+1}^+)\in V_h|_{\Omega_h[\phi_h(t_{m+1}^-)]}$ is computed from
\eqref{eq:scheme-u}.

\section{Main result}
\label{sec:main}

\begin{assumption}[on the exact solution]
\label{ass:exact}
Let $n\in\{2,3\}$ and $k\ge3$. The extended bulk solution $u$, the source $f$, the level-set function $\phi$, and the response function $\beta$ satisfy \eqref{eq:bulk} and \eqref{eq:levelset-domain}--\eqref{eq:HJ}, and
\begin{equation}
\begin{aligned}
 \phi &\in C^2([0,T];W^{k+2,\infty}(D)),\qquad
 &&u\in C^2([0,T];W^{k+2,\infty}(D)),\\
 f&\in C^1([0,T];W^{k,\infty}(D)),\qquad
 &&\beta\in C^{k+2}(\R).
\end{aligned}
 \label{eq:standing-smoothness}
\end{equation}
\end{assumption}

We split the errors into a discrete part and an interpolation part,
\begin{equation}
 e_w:=w_h-I_hw \quad\mbox{and}\quad \rho_w:=I_hw-w\quad
 \mbox{for}\,\,\,w\in\{u,p,\eta,\phi\},
 \label{eq:discrete-errors}
\end{equation}
so that $w_h-w=e_w+\rho_w$ and $\rho_p=\rho_\phi$. Thus $e_u\in V_h|_{\Oha}$, $e_p,e_\phi\in V_h|_{\Uha}$ and $e_\eta\in(V_h|_{\Uha})^n$.

\begin{theorem}[Error estimate]
\label{thm:main}
Under \cref{ass:exact}, there exists $h_0>0$ such that for $0<h\le h_0$ the scheme \eqref{eq:scheme-p}--\eqref{eq:transfer} has a solution $(p_h,\eta_h,\phi_h,u_h)$, unique on the local branch, which is continuous on each interval $I_m$, and the update times satisfy \eqref{eq:event-time interval-spacing}. Moreover,
\begin{equation}
 \sup_{0\le t\le T}
 \Bigl(
 \nH{e_u}{\Oha} +  \nH{e_\phi}{\Uha}
 + \nL{e_p}{\Uha} +\nL{e_\eta}{\Uha}
 \Bigr)
 \le Ch^{k}.
 \label{eq:main-theorem-error}
\end{equation}
At an update time $t_m$, the supremum includes both $t_m^-$ and $t_m^+$, each measured on its own active meshes.
\end{theorem}

\begin{remark}
\label{rem:main-consequences}
By \eqref{eq:main-theorem-error}, the discrete Sobolev inequality \eqref{eq:discrete-sobolev} and \eqref{eq:Nd-def}, the Hausdorff distance between $\Gamma[\phi_h(t)]$ and $\Gamma(t)$ is at most $Ch^{k-(n-1)/4}$.
\end{remark}

\begin{remark}[Error estimate for the fully discrete scheme]
\label{rem:fd-estimate}
We expect that the analysis extends to the fully discrete semi-implicit scheme of \cref{subsec:fd-scheme} under \cref{ass:exact} and the step-size condition $\tau\le ch^k$. We do not pursue this here and prove \cref{thm:main} only for the semidiscrete scheme.
\end{remark}

\section{Preliminaries}
\label{sec:prelim}

Throughout the analysis, $(\cdot,\cdot)_S$ denotes the $L^2$ inner product on a domain or hypersurface $S$. For $q\in\R^n\setminus\{0\}$, set
\begin{equation}
 n(q):=\frac q{|q|},\qquad P(q):=I-n(q)\otimes n(q),\qquad
 Dn(q)=\frac1{|q|}P(q).
 \label{eq:Dn}
\end{equation}
In particular, for $q=\nabla\phi$, the quantities in \eqref{eq:Dn} are the normal vector, the tangential projection, and the derivative of the normalization map $n:\R^n\setminus\{0\}\rightarrow \R^n$. 

Unless stated otherwise, we denote by $C$ a generic constant that is independent of $h$, the cut position, and the update index.  It may depend on the exact solution, $T$, the fixed parameters of the method, and the uniform mesh constants. 

\subsection{Interpolation and inverse estimates}
\label{subsec:interp-geometry}

We collect the interpolation and inverse estimates needed for the error analysis. 
For $w\in W^{k+1,\infty}(D)$ the nodal interpolant satisfies the standard estimates (see \cite[Theorem~3.1.6]{Ciarlet1978} or 
\cite[Theorem~4.4.20]{BrennerScott2008}): 
\begin{equation}
 \|w-I_hw\|_{L^\infty(D)}+h\|w-I_hw\|_{W^{1,\infty}(D)}\le Ch^{k+1}\|w\|_{W^{k+1,\infty}(D)} . 
 \label{eq:phi-p-interpolation-Winf}
\end{equation}
Consequently, under \Cref{ass:exact}, the interpolation errors \eqref{eq:discrete-errors} satisfy, uniformly for $t\in[0,T]$,
\begin{equation}
\begin{aligned}
 \|\rho_\phi\|_{L^\infty(D)}+\|\rho_p\|_{L^\infty(D)}+\|\rho_u\|_{L^\infty(D)}
 +\|\rho_\eta\|_{L^\infty(D)} &\le Ch^{k+1},\\
 \|\rho_\phi\|_{W^{1,\infty}(D)}+\|\rho_u\|_{W^{1,\infty}(D)}
 +\|\rho_\eta\|_{W^{1,\infty}(D)} &\le Ch^{k},\\
 \|\partial_t\rho_\phi\|_{L^\infty(D)}+\|\partial_t\rho_\eta\|_{L^\infty(D)}
 &\le Ch^{k+1} .
\end{aligned}
 \label{eq:interp-ass}
\end{equation}

The estimate \eqref{eq:interp-ass} guarantees that the interpolated level set $\phi_h^\ast:=I_h\phi$ inherits the nondegeneracy of $\phi$ for $h$ small: $\tfrac12c\le|\nabla\phi_h^\ast|\le C$ elementwise in $\Nd(t)$.  For $c\in\{-\delta,0,\kappa h\}$ and sufficiently small $h$, its level surfaces satisfy $\dist(x,\Gamma_c[\phi])\le Ch^{k+1}$ for $x\in \Gamma_c[\phi_h^\ast]$. 

For finite element functions we use the inverse inequalities (see \cite[Theorem~3.2.6]{Ciarlet1978} or \cite[Theorem~4.5.11]{BrennerScott2008}): for $v_h\in V_h$, $K\in\T_h$ and $1\le q\le\infty$,
\begin{equation}
\begin{gathered}
 \nL{\nabla v_h}{K}\le Ch^{-1}\nL{v_h}{K},\qquad
 \|v_h\|_{L^q(K)}\le Ch^{-n(\frac12-\frac1q)}\nL{v_h}{K},\\
 \nLi{\nabla v_h}{K}\le Ch^{-1-\frac n2}\nL{v_h}{K}.
\end{gathered}
 \label{eq:inverse-inequalities}
\end{equation}

\subsection{Level-set shape calculus}
\label{subsec:shape}

We compare integrals on bands defined by level-set functions $\psi^\theta$, $\theta\in[0,1]$, whose relevant level surfaces are uniformly nondegenerate Lipschitz graphs, compactly contained in the fixed active mesh $\Uha$, and depend Lipschitz continuously on $\theta$. Applying the coarea formula and the fundamental theorem of calculus elementwise gives the following integrated identity (see also \cite[Chapter~9]{DelfourZolesio2011}):
\begin{align}\label{eq:FTC-band}
 \int_{\Ub{\psi^1}}F^1\,\dd x-\int_{\Ub{\psi^0}}F^0\,\dd x
 =\int_0^1\Bigl[\int_{\Ub{\psi^\theta}}\partial_\theta F^\theta\,\dd x
 &-\int_{\Gamma_+[\psi^\theta]}F^\theta\,\frac{\partial_\theta\psi^\theta}{|\nabla\psi^\theta|}\,\dd s \notag\\
 &+\int_{\Gamma_-[\psi^\theta]}F^\theta\,\frac{\partial_\theta\psi^\theta}{|\nabla\psi^\theta|}\,\dd s\Bigr]\dd\theta .                                                
\end{align}
It is the analogue of the Reynolds transport theorem for the band $\Ub{\phi(t)}$:
\begin{equation}
 \frac{\dd}{\dd t}\int_{\Ub{\phi(t)}}F\,\dd x
 =\int_{\Ub{\phi(t)}}\partial_tF\,\dd x
 -\int_{\Gamma_+[\phi]}F\frac{\partial_t\phi}{|\nabla\phi|}\,\dd s
 +\int_{\Gamma_-[\phi]}F\frac{\partial_t\phi}{|\nabla\phi|}\,\dd s ,
 \label{eq:reynolds-band}
\end{equation}
where $\Gamma_\pm[\phi]$ move with outward normal velocity $\mp\partial_t\phi/|\nabla\phi|=\pm\beta(u)$.

For the bulk error analysis we use the level-set homotopy
\begin{equation}
 \phi_h^\theta:=\phi_h^\ast+\theta e_\phi
 \label{eq:state-theta}
\end{equation}
with $\phi_h^\ast:=I_h \phi$, so that 
$$\phi_h^0=\phi_h^\ast,\quad \phi_h^1=\phi_h,\quad 
\partial_\theta\phi_h^\theta=e_\phi.$$

\subsection{The continuation interval and associated bounds}
\label{subsec:continuation}

We define $t^\ast$ to be the supremum of the times $t\in(0,T]$ for which the semidiscrete solution on the local branch exists on $[0,t)$ and
satisfies, for every $s\in[0,t)$,
\begin{equation}
\begin{aligned}
 \nH{e_u(s)}{\Oha}+\nH{e_\phi(s)}{\Uha}
 &\le h^{k-0.1},\\
 \nL{e_p(s)}{\Uha}+\nL{e_\eta(s)}{\Uha}
 &\le h^{k-0.1}.
\end{aligned}
 \label{eq:continuation-conditions}
\end{equation}
We set the supremum to zero if there is no such positive candidate time. At each update time $t_{m+1}$ strictly smaller than the candidate time $t$, these inequalities are required for both one-sided values at $t_{m+1}^\pm$, with the errors measured on their respective active meshes. All subsequent a priori error estimates are understood on $[0,t^\ast)$, with constants independent of $t^\ast$; the positivity of $t^\ast$ and the validity of \eqref{eq:continuation-conditions} at the left limit at $t^\ast$ are justified in \cref{sec:closure}. 

In the following, the statement ``\emph{sufficiently small $h$}'' means smaller than some constant which is independent of $t^*$ (possibly depending on $T$). The inverse estimates \eqref{eq:inverse-inequalities} imply that, for $t\in[0,t^\ast)$ and sufficiently small $h$, 
\begin{align}
 \nWi{e_\phi}{\Uha}
 &\le C h^{-\frac n2}\nH{e_\phi}{\Uha}
 \le Ch^{k-0.1-\frac n2} \leq h  \label{eq:ephi-Winf-from-H1} ,
\end{align}
where the last inequality uses $k\ge3$ and $n\le3$. The same elementwise inverse estimate gives
$$
 \nWi{e_u}{\Oha}\le Ch^{-n/2}\nH{e_u}{\Oha}
 \le Ch^{k-0.1-n/2},
$$
and
\begin{equation}
 \nLi{e_p}{\Uha}+\nLi{e_\eta}{\Uha}
 \le Ch^{k-0.1-n/2}.
 \label{eq:epsi-Linf}
\end{equation}
Since $k-0.1-n/2\ge1.4$, these ordinary inverse estimates already imply, independently of \cref{lem:discrete-sobolev},
\begin{align}
\begin{aligned}
 \nLi{e_u}{\Oha} + \nLi{e_\phi}{\Uha} &\le h^{1.1}, \\
 \nLi{e_p}{\Uha}+\nLi{e_\eta}{\Uha} &\le h^{1.1} .
\end{aligned}
 \label{eq:smallness}
\end{align}
These estimates imply $\nWi{u_h}{\Oha}\le C$ and $\nLi{\phi_h-\phi}{\Uha}\le Ch^{1.1}$. 

By the upper-distance condition in \eqref{eq:band-buffer}, any $x\in\Uha$ is within $2c_{\rm buf}h$ of a point $y\in\overline{\Uh}$. The global regularity of $\phi$ gives
$$
 |\phi(x)|\le\max\{\delta,\kappa h\}
       +\nLi{\phi_h-\phi}{\Uha}+Ch<2\delta
$$
for small $h$. Thus $\Uha\subset\Nd(t)$. For sufficiently small $h$, \eqref{eq:ephi-Winf-from-H1} and \eqref{eq:smallness} ensure that each $\Gamma_c[\phi_h(t)]$ with $c\in\{-\delta,0,\kappa h\}$  is a small Lipschitz normal graph over $\Gamma_c[\phi(t)]$, and
\begin{align}\label{error-distance}
\begin{aligned}
 \dist(x,\Ga{\phi})&\le Ch^{1.1}
       &&\text{for }x\in\Gh,\\
 \dist(x,\Gamma_\pm[\phi])&\le Ch^{1.1}
       &&\text{for }x\in\Gamma_\pm[\phi_h],\\
 \dist(x,\Om{\phi})&\le Ch^{1.1}
       &&\text{for }x\in\Oh.
\end{aligned}
\end{align}
The first two inequalities also hold with the numerical and exact surfaces interchanged. The sign correspondence and the normal segments give
\begin{equation}
 \Ga{\phi}\subset\{|\phi_h|<\tfrac12\kappa h\},\qquad
 \Om{\phi}\subset\Oha,\qquad
 \Ub{\phi}\subset\Uha\subset\Nd(t).
 \label{eq:all-bands-in-active}
\end{equation}

\subsection{The chain condition}
\label{section:patch}

For $t\in I_m\cap[0,t^\ast)$ and sufficiently small $h$, the definitions in \eqref{eq:buffer-elements}--\eqref{eq:numerical-bulk-support}
ensure that every element $K$ in $\Uha$ contains a point in the $c_{\rm buf}h$-neighborhood of $U[\phi_h(t_m^-)]$. Joining this point to a closest point in $\overline{U[\phi_h(t_m^-)]}$, and preceding this segment by a segment inside $K$ if necessary, gives a path from the interior of $K$ to the reference band $U[\phi_h(t_m^-)]$ of length at most $Ch$ contained in $\Uha$. This gives the chain condition of \cref{def:chain}; see \cref{appendix:path}. 

\subsection{Lower bound of \texorpdfstring{$t_{m+1}-t_m$ (time between active-mesh updates)}{the time between active-mesh updates}}
\label{subsec:interval}
According to the definitions in \eqref{eq:buffer-elements}--\eqref{eq:numerical-bulk-support}, the active mesh $U_h[\phi_h(t_m^-)]$ contains the $c_{\rm buf}h$-neighborhood of $U[\phi_h(t_m^-)]$ and is contained in its $(c_{\rm buf}+1)h$-neighborhood. Therefore, there is a gap of $\frac12c_{\rm buf}h$ between the $\frac12c_{\rm buf}h$-neighborhood of $U[\phi_h(t_m^-)]$ and the boundaries of $U_h[\phi_h(t_m^-)]$. Then the bounded speed of the moving band $U[\phi]$ and the $O(h^{1.1})$ distance between $U[\phi]$ and $U[\phi_h]$ ensure that \eqref{eq:band-buffer-a} cannot fail within $c_{\rm interval}^*h$ time for some small positive constant $c_{\rm interval}^*$. Similarly, there is a $(c_{\rm buf}-1)h$-gap between $U_h[\phi_h(t_m^-)]$ and the complement of the $2c_{\rm buf}h$-neighborhood of $U[\phi_h(t_m^-)]$. Hence \eqref{eq:band-buffer-b} cannot fail within $c_{\rm interval}^\# h$ time for some small positive constant $c_{\rm interval}^\#$. This means that, for $c_{\rm interval}=\min(c_{\rm interval}^* ,c_{\rm interval}^\# )$, 
\begin{equation}\label{time-length-min}
 t_{m+1}-t_m\ge c_{\rm interval}h
 \qquad\text{for every update }t_{m+1}<\min\{t^\ast,T\}.
\end{equation}
If $t^*<T$ coincides with an update time $t_{m+1}$ then the same conclusion holds. Hence there are at most $C/h$ active-mesh updates. 

In addition to the inverse estimates in \eqref{eq:inverse-inequalities}, we use the following bound of the $L^\infty$ norm of a finite element function on $U_h[\phi_h]$ by its $H^1$ norm on $U_h[\phi_h]$.
\begin{lemma}[Discrete Sobolev inequality]
\label{lem:discrete-sobolev}
Under \cref{ass:exact}, the following estimate holds for sufficiently small $h$ and $t\in[0,t^\ast)$: 
\begin{equation}
 \|v_h\|_{L^\infty(U_h[\phi_h])}
 \le
 Ch^{-(n-1)/4}\|v_h\|_{H^1(U_h[\phi_h])}
 \quad
 \forall v_h\in V_h|_{U_h[\phi_h]},
 \quad n\in\{2,3\},
 \label{eq:discrete-sobolev}
\end{equation}
where $C$ is independent of $h$, the update index, and the cut position.
\end{lemma}

\begin{remark}
In two dimensions the classical discrete Sobolev inequality gives a logarithmic factor rather than $h^{-1/4}$; however, the present bound suffices for the following proof. 
\end{remark}

\begin{proof}
By the estimates \eqref{eq:interp-ass} and \eqref{eq:ephi-Winf-from-H1}, 
$
 \|\phi_h-\phi\|_{W^{1,\infty}(U_h[\phi_h])}
 \le Ch^{k-0.1-\frac{n}{2}}.
$
For $h$ sufficiently small, the two boundaries of $U[\phi_h]$ are uniformly Lipschitz graphs over the corresponding boundaries of $U[\phi]$. The fixed interior width of the band therefore ensures a uniformly bounded extension operator $E_{U,D}$ from $H^1(U[\phi_h])$ to $H^1(D)$. Consequently, the continuous Sobolev embedding yields, for $q=8$ in two space dimensions and $q=6$ in three, 
\begin{align}
\begin{aligned}
 \|v_h\|_{L^q(U[\phi_h])}
 \leq \|E_{U,D}v_h\|_{L^q(D)}
 \leq C\|E_{U,D}v_h\|_{H^1(D)}
 &\le C\|v_h\|_{H^1(U[\phi_h])}\\
 &\le C\|v_h\|_{H^1(U_h[\phi_h])} .
\end{aligned}
\label{eq:band-sobolev-Lq}
\end{align}

We next control the values of $v_h$ in the buffer region $U_h[\phi_h]\backslash U[\phi_h]$. By the path construction in \cref{appendix:path}, every point $x\in U_h[\phi_h]$ can be joined to a ball $Q_x\subset U[\phi_h]$, with $|Q_x|\simeq h^n$, by paths of length at most $Ch$ contained in $\overline{U_h[\phi_h]}$. 

Since $v_h$ is continuous and piecewise polynomial, integration along these paths gives
$$
 |v_h(x)|
 \le |v_h(y)|
      +Ch\|\nabla v_h\|_{L^\infty(U_h[\phi_h])}
 \qquad \forall y\in Q_x.
$$
Averaging over $Q_x$ and applying H\"older's inequality and the elementwise inverse estimate, we obtain
\begin{align*}
 |v_h(x)|
 &\le
 Ch^{-n/q}\|v_h\|_{L^q(Q_x)}
 +Ch\|\nabla v_h\|_{L^\infty(U_h[\phi_h])}\\
 &\le
 Ch^{-n/q}\|v_h\|_{L^q(U[\phi_h])}
 +Ch^{1-n/2}\|\nabla v_h\|_{L^2(U_h[\phi_h])}.
\end{align*}
Taking the supremum over $x$ and using \eqref{eq:band-sobolev-Lq} gives
$$
 \|v_h\|_{L^\infty(U_h[\phi_h])}
 \le
 C\bigl(h^{-n/q}+h^{1-n/2}\bigr)
 \|v_h\|_{H^1(U_h[\phi_h])}.
$$
For the stated choices of $q$ and $0<h\le1$, we have $ h^{-n/q}+h^{1-n/2}
 \le Ch^{-(n-1)/4}$, which proves \eqref{eq:discrete-sobolev}.
\end{proof}

We may now sharpen the preliminary value bound by applying
\eqref{eq:discrete-sobolev} to $e_\phi$:
\begin{equation}
 \nLi{e_\phi}{\Uha}
 \le Ch^{-(n-1)/4}\nH{e_\phi}{\Uha}
 \le Ch^{k-0.1-(n-1)/4}.
 \label{eq:ephi-Linf-W14}
\end{equation}
The preceding localization and path arguments did not use this sharper
estimate.

\section{Error equations}
\label{sec:error-eq}

Throughout this section and the next, $t\in[0,t^\ast)$ and $h$ is
small enough for \eqref{eq:ephi-Winf-from-H1}--\eqref{eq:discrete-sobolev}.
Estimates involving time derivatives are used only inside the intervals $I_m$, where the active meshes are fixed; the other estimates also hold for the one-sided states at update times.

\subsection{The exact fields on the numerical geometry}
\label{subsec:exact-identities}

Since a globally $C^k$ function has no jumps of its normal derivatives across element faces,
\begin{equation}
 g_r(w,v_h)_{\Uha}=0,\qquad r=0,1,
 \label{eq:ghost-smooth-zero}
\end{equation}
for every $w\in C^k(\bar D)$, in particular for $u$, $p$, $\eta$, $\phi$ and their time derivatives.  The standard face-jump interpolation estimate with the same scalings also yields, uniformly on the active meshes, $g_1(I_hw,I_hw)^{1/2}\le Ch^k\|w\|_{W^{k+1,\infty}(D)}$ and $g_0(I_hw,I_hw)^{1/2}\le Ch^{k+1}\|w\|_{W^{k+1,\infty}(D)}$; these bounds will be used without further comment. 

Integrating the extended identity $u-\Delta u=f + r$, where
$r=u-\Delta u-f$, against $v_h\in V_h|_{\Oha}$ on $\Oh$ and
integrating by parts, we obtain
\begin{equation}
\begin{aligned}
 &(u,v_h)_{\Oh}+(\nabla u,\nabla v_h)_{\Oh}+g_1(u,v_h)_{\Uha}\\
 &\quad=(f(t),v_h)_{\Oh}+(r,v_h)_{\Oh}+(n(\nabla\phi_h)\cdot\nabla u,v_h)_{\Gh},
\end{aligned}
 \label{eq:exact-u-on-oh}
\end{equation}
where $n(\nabla\phi_h)$ is the outward unit normal vector on $\Gh=\partial \Oh$. The last two terms are the geometric consistency error; they arise because \eqref{eq:bulk} holds on $\Om{\phi}$, not on $\Oh$. Since \eqref{eq:p-cont} and \eqref{eq:eta-cont} hold pointwise in $\Nd(t)\supset\Uha$, the following equations hold for all $v_h\in V_h|_{\Uha}$:
\begin{align}
 &(\partial_t p+\beta(u_h)n(\eta_h)\cdot\nabla p,v_h)_{\Uh}
 +g_0(\partial_tp+p,v_h)_{\Uha}\nonumber\\
 &\quad=\bigl((\beta(u_h)n(\eta_h)-\beta(u)n(\eta))\cdot\nabla p,v_h\bigr)_{\Uh},
 \label{eq:exact-p-on-uh}\\
 &(\partial_t\eta_j+\beta(u_h)n(\eta_h)\cdot\nabla\eta_j,v_h)_{\Uh}
 +g_0(\partial_t\eta_j+\eta_j,v_h)_{\Uha}\nonumber\\
 &\quad=-(\beta'(u)\partial_ju|\nabla\phi|,v_h)_{\Uh}
 +\bigl((\beta(u_h)n(\eta_h)-\beta(u)n(\eta))\cdot\nabla\eta_j,v_h\bigr)_{\Uh} . 
 \label{eq:exact-eta-on-uh}
\end{align}
Moreover, since $p=\phi$ and $\eta=\nabla\phi$, it follows that 
\begin{equation}
\begin{aligned}
 &(\phi,v_h)_{\Uh}+(\nabla\phi,\nabla v_h)_{\Uh}+g_1(\phi,v_h)_{\Uha}\\
 &\qquad=(p,v_h)_{\Uh}+(\eta,\nabla v_h)_{\Uh}
 \qquad\forall v_h\in V_h|_{\Uha} .
\end{aligned}
 \label{eq:exact-phi-on-uh}
\end{equation}

\subsection{The bulk error identity}
\label{subsec:bulk-error-identity}

Subtract \eqref{eq:exact-u-on-oh} from \eqref{eq:scheme-u},
use $u_h=u+\rho_u+e_u$, and take $v_h=e_u$. Since the same band-face
penalty is used in both equations, this gives
\begin{align} \label{eq:error-u-two-geometry}
\begin{aligned}
 &\|e_u\|_{L^2(\Omega[\phi_h])}^2+
   \|\nabla e_u\|_{L^2(\Omega[\phi_h])}^2 + g_1(e_u,e_u)_{\Uha}\\
 &\quad=-(r,e_u)_{\Omega[\phi_h]}
 -(n(\nabla\phi_h)\cdot\nabla u,e_u)_{\Gamma[\phi_h]}\\
 &\qquad-(\rho_u,e_u)_{\Omega[\phi_h]}
 -(\nabla\rho_u,\nabla e_u)_{\Omega[\phi_h]}
 -g_1(\rho_u,e_u)_{\Uha}.
\end{aligned}
\end{align}
The residual $r=u-\Delta u-f\in W^{k,\infty}(D)$ vanishes in
$\Omega[\phi]$, so its derivatives through order $k-1$ vanish on
$\Gamma[\phi]$. Taylor's formula along an exterior normal segment gives
$|r(x)|\le C\dist(x,\Gamma[\phi])^k$, while $r=0$ inside.
Together with \eqref{error-distance}, this yields
\begin{equation}
|r|\leq Ch^k \quad\mbox{on}\,\,\,\,
\Omega[\phi_h] .
\label{eq:r-shape-difference}
\end{equation}

\subsection{The transport error equations}
\label{subsec:transport-error-eq}

Subtracting \eqref{eq:exact-p-on-uh} from \eqref{eq:scheme-p} and using $p_h=p+\rho_p+e_p$, we obtain the following error equation for $e_p$:
\begin{equation}
\begin{aligned}
 &(\partial_te_p+\beta(u_h)n(\eta_h)\cdot\nabla e_p,v_h)_{\Uh}
 +g_0(\partial_te_p+e_p,v_h)_{\Uha}\\
 &\quad=-(\partial_t\rho_p+\beta(u_h)n(\eta_h)\cdot\nabla\rho_p,v_h)_{\Uh}\\
 &\qquad-\bigl((\beta(u_h)n(\eta_h)-\beta(u)n(\eta))\cdot\nabla p,v_h\bigr)_{\Uh}\\
 &\qquad-g_0(\partial_t\rho_p+\rho_p,v_h)_{\Uha} \qquad\forall\, v_h\in V_h|_{\Uha} .
\end{aligned}
 \label{eq:error-p-current-band}
\end{equation}
Similarly, subtracting \eqref{eq:exact-eta-on-uh} from \eqref{eq:scheme-eta} and using $\eta_h=\eta+\rho_\eta+e_\eta$, we obtain the following error equation for $e_\eta$: 
\begin{equation}
\begin{aligned}
 &(\partial_te_{\eta}+\beta(u_h)n(\eta_h)\cdot\nabla e_{\eta},\zeta_h)_{\Uh}
 +g_0(\partial_te_{\eta}+e_{\eta},\zeta_h)_{\Uha}\\
 &\quad=-(\partial_t\rho_{\eta}+\beta(u_h)n(\eta_h)\cdot\nabla\rho_{\eta},\zeta_h)_{\Uh}\\
 &\qquad-\bigl((\beta(u_h)n(\eta_h)-\beta(u)n(\eta))\cdot\nabla\eta,\zeta_h\bigr)_{\Uh}\\
 &\qquad-\bigl(\beta'(u_h)\nabla u_h|\eta_h|-\beta'(u)\nabla u|\eta|,\zeta_h\bigr)_{\Uh}\\
 &\qquad-g_0(\partial_t\rho_{\eta}+\rho_{\eta},\zeta_h)_{\Uha} \qquad\forall\, \zeta_h\in (V_h|_{\Uha})^n .
\end{aligned}
 \label{eq:eeta-current-band}
\end{equation}

Since \eqref{eq:smallness} implies
$$
 \|\beta(u_h)n(\eta_h)\|_{L^\infty(\Uha)}\le C,
$$
\eqref{eq:interp-ass} gives
\begin{align}
& \nLi{\partial_t\rho_p+\beta(u_h)n(\eta_h)\cdot\nabla\rho_p}{\Uha} \notag\\
&\quad +\nLi{\partial_t\rho_\eta+\beta(u_h)n(\eta_h)\cdot\nabla\rho_\eta}{\Uha} \le Ch^k . 
 \label{eq:Dm-rho-bound}
\end{align}
The face-jump interpolation bound following \eqref{eq:ghost-smooth-zero}, applied to $p$, $\eta$ and their time derivatives, gives  
\begin{equation}
 g_0(\partial_t\rho_p+\rho_p,\partial_t\rho_p+\rho_p)_{\Uha}
 + g_0(\partial_t\rho_\eta+\rho_\eta,\partial_t\rho_\eta+\rho_\eta)_{\Uha}
 \le Ch^{2k+2}.
 \label{eq:ghost-rho-bound}
\end{equation}

\subsection{The reconstruction error equation}
\label{subsec:rec-error-eq}

Subtracting \eqref{eq:exact-phi-on-uh} from \eqref{eq:scheme-rec} and using $\phi_h=\phi+\rho_\phi+e_\phi$, $p_h=p+\rho_p+e_p$, $\eta_h=\eta+\rho_\eta+e_\eta$ and $\rho_p=\rho_\phi$, we obtain the following error equation for $e_\phi$: 
\begin{equation}
\begin{aligned}
 &(e_\phi,v_h)_{\Uh}+(\nabla e_\phi,\nabla v_h)_{\Uh}+g_1(e_\phi,v_h)_{\Uha}\\
 &\quad=(e_p,v_h)_{\Uh}+(e_\eta,\nabla v_h)_{\Uh}
 +(\rho_\eta-\nabla\rho_\phi,\nabla v_h)_{\Uh}
 -g_1(\rho_\phi,v_h)_{\Uha}, 
\end{aligned}
 \label{eq:ephi-current-band}
\end{equation}
which holds for all $v_h\in V_h|_{\Uha}$. 
By \eqref{eq:interp-ass}, \eqref{eq:ghost-smooth-zero} and the face-jump interpolation bound following it,
\begin{equation}
 \nLi{\rho_\eta-\nabla\rho_\phi}{D}\le Ch^k,\quad
 g_1(\rho_\phi,\rho_\phi)_{\Uha}=g_1(I_h\phi,I_h\phi)_{\Uha}\le Ch^{2k}.
 \label{eq:ephi-interp-residual}
\end{equation}

\section{Error estimates for \texorpdfstring{$t\in[0,t^\ast)$}{the continuation interval}}
\label{sec:apriori}

Starting from the error identities in \cref{sec:error-eq}, we
derive the estimates used in \cref{sec:closure} to improve the continuation bounds in \eqref{eq:continuation-conditions}. 

\subsection{Estimates for \texorpdfstring{${\bf e_u}$}{the bulk error}}
\label{subsec:eu-estimate}

Let $\Gamma:=\Ga{\phi}$ and let $\pi_h\colon\Gh\to\Gamma$ be the closest-point map, which is a bi-Lipschitz bijection for sufficiently small $h$ by the normal-graph property in \cref{subsec:continuation}, together with \eqref{eq:interp-ass} and \eqref{eq:ephi-Winf-from-H1}. For a function $v$ on $\Gh$ we denote by $v^\ell:=v\circ\pi_h^{-1}$ its lift to $\Gamma$.

\begin{lemma}[Boundary and residual terms]
\label{lem:bulk-terms}
The following estimates hold for sufficiently small $h$ and $t\in[0,t^\ast)$. 
\begin{enumerate}
\item[{\rm(}a{\rm)}] {\rm(}Exact-normal derivative on the numerical interface{\rm)} 
\begin{equation}
 |(n(\nabla\phi)\cdot\nabla u,e_u)_{\Gh}|
 \le C\bigl(\nH{e_\phi}{\Uha}+h^{k+1}\bigr)\nH{e_u}{\Oha} .                   \label{eq:exact-normal-bound}
\end{equation}
\item[{\rm(}b{\rm)}] {\rm(}Error between normal vectors{\rm)} The following relations hold on $\Gh${\rm:}
\begin{align}
 &n(\nabla\phi_h)-n(\nabla\phi) =\int_0^1\frac{P(\nabla\phi_h^\theta)}{|\nabla\phi_h^\theta|}\nabla e_\phi\,\dd\theta
 +\bigl[n(\nabla\phi_h^\ast)-n(\nabla\phi)\bigr], \label{eq:normal-FTC}\\
 &\bigl((n(\nabla\phi_h)-n(\nabla\phi))\cdot\nabla u,e_u\bigr)_{\Gh} =T_h+R_h, \label{eq:normal-TR}
\end{align} 
with $T_h:=(|\nabla\phi|^{-1}P(\nabla\phi_h)\nabla e_\phi\cdot\nabla u,e_u)_{\Gh}$ and 
\begin{align}
|R_h| \le Ch^{\frac12}\| e_\phi\|_{H^1(U_h[\phi_h])} \|e_u\|_{H^1(\Omega_h[\phi_h])} + Ch^k \|e_u\|_{H^1(\Omega_h[\phi_h])} .
 \label{eq:normal-remainder-bound}
\end{align}
\item[{\rm(}c{\rm)}] {\rm(}Lift and duality{\rm)} 
Moreover, $T_h=(|\nabla\phi|^{-1}\nabla_\Gamma u\cdot\nabla_\Gamma e_\phi^\ell,e_u^\ell)_\Gamma+R_{\rm lift}$, with
\begin{equation}
\begin{aligned}
 |R_{\rm lift}|&\le C\bigl(\nWi{e_\phi}{\Uha}+h^k\bigr)h^{-\frac12}\nH{e_\phi}{\Uha}\nH{e_u}{\Oha} 
\end{aligned}                                                        \label{eq:lift-remainder}
\end{equation}
and
\begin{equation}
 \bigl|\bigl(|\nabla\phi|^{-1}\nabla_{\Gamma[\phi]} u\cdot\nabla_{\Gamma[\phi]} e_\phi^\ell,e_u^\ell\bigr)_{\Gamma[\phi]}\bigr|
 \le C\nH{e_\phi}{\Uha}\nH{e_u}{\Oha} .                            \label{eq:tangential-main-bound}
\end{equation}
\end{enumerate}
\end{lemma}

\begin{proof}
(a) The function $q:=n(\nabla\phi)\cdot\nabla u$ is zero on $\Gamma$ because of the Neumann condition on $u$, and therefore $|q(x)|\le C\dist(x,\Gamma)\le C|\phi(x)|$ on $\Nd(t)$. Since $\phi_h=0$ on $\Gh$, substituting $\phi=-\rho_\phi-e_\phi$ yields $|q|\le C(|e_\phi|+|\rho_\phi|)$ on $\Gh$. The trace estimates follow by restricting both errors to $U[\phi_h]\cap\Omega[\phi_h]$, whose interior width is fixed and whose boundary contains $\Gh$. Its uniform Lipschitz charts bound the $L^2$ and $H^{1/2}$ traces on $\Gh$ by the respective active-mesh $H^1$ norms, independently of the cut position. This proves (a) and supplies the traces used below. 

(b) Identity \eqref{eq:normal-FTC} is the fundamental theorem of calculus for $\theta\mapsto n(\nabla\phi_h^\theta)$ with \eqref{eq:Dn}. The remainder 
$
R_h
$
can be written $R_h = R_h^{(1)}+R_h^{(2)}$ with
$$
\begin{aligned}
R_h^{(1)}&:=\Big(\Big[\int_0^1\frac{P(\nabla\phi_h^\theta)}{|\nabla\phi_h^\theta|}\dd\theta-\frac{P(\nabla\phi_h)}{|\nabla\phi|}\Big]\nabla e_\phi\cdot\nabla u,\ e_u\Big)_{\Gamma[\phi_h]},\\
R_h^{(2)}&:=\Big(\big[n(\nabla\phi_h^*)-n(\nabla\phi)\big]\cdot\nabla u,\ e_u\Big)_{\Gamma[\phi_h]} .
\end{aligned}
$$
Since $|n(\nabla\phi_h^\ast)-n(\nabla\phi)|\leq Ch^k$ and $\big|\frac{P(\nabla\phi_h^\theta)}{|\nabla\phi_h^\theta|}-\frac{P(\nabla\phi_h)}{|\nabla\phi|}\big|\leq C(|\nabla e_\phi|+|\nabla \rho_\phi|)$, 
we obtain
\begin{align*}
|R_h| \le &C\int_{\Gh}(|\nabla e_\phi|^2+|\nabla\rho_\phi||\nabla e_\phi|)|e_u|\dd s + Ch^k\|e_u\|_{L^2(\Gh)} \\
&\quad \leq Ch\|\nabla e_\phi\|_{L^2(\Gamma[\phi_h])} \|e_u\|_{L^2(\Gamma[\phi_h])} + Ch^k\|e_u\|_{L^2(\Gh)} 
\quad\mbox{(here \eqref{eq:ephi-Winf-from-H1} is used)}\\
&\quad \leq Ch^{\frac12}\| e_\phi\|_{H^1(U_h[\phi_h])} \|e_u\|_{H^1(\Omega_h[\phi_h])} + Ch^k \|e_u\|_{H^1(\Omega_h[\phi_h])} , 
\end{align*}
where the last inequality uses the inverse trace estimate for $\nabla e_\phi$ and the trace bound from {\rm(}a{\rm)} for $e_u$. 

(c) With $\bm w:=|\nabla\phi|^{-1}\nabla u$ and $\nabla_{\Gamma[\phi_h]}=P(\nabla\phi_h)\nabla$ on $\Gamma[\phi_h]$,
$$
R_{\rm lift}=(\bm w\cdot\nabla_{\Gamma[\phi_h]}e_\phi,e_u)_{\Gamma[\phi_h]}-(\bm w\cdot\nabla_\Gamma e_\phi^\ell,e_u^\ell)_\Gamma .
$$
Transform the second term to $\Gamma[\phi_h]$ with $y=\pi_h(x)$. Using $(\nabla_\Gamma e_\phi^\ell)\circ\pi_h=B_h\nabla_{\Gamma[\phi_h]}e_\phi$, with $B_h:=(P(\nabla\phi_h)D\pi_h^T)^{-1}$ as a map from the tangent plane at $x\in\Gh$ to the tangent plane at $\pi_h(x)\in\Gamma$, it becomes $\int_{\Gamma[\phi_h]}(\bm w\circ\pi_h)\cdot(B_h\nabla_{\Gamma[\phi_h]}e_\phi)\,e_u\,J_h\dd s$, where $J_h$ is the surface Jacobian of $\pi_h$. In the operator norm below, $B_h$ and $P(\nabla\phi_h)$ are restricted to the tangent plane at $x\in\Gh$; no inverse of an ambient singular matrix is used. On $\Gamma[\phi_h]$ we have $|\bm w-\bm w\circ\pi_h|\le C\dist(\cdot,\Gamma)$ and
$$
\begin{aligned}
|J_h-1|+|B_h-P(\nabla\phi_h)|&\le C\big(\dist(\cdot,\Gamma)+|n(\nabla\phi_h)-n(\nabla\phi)\circ\pi_h|\big)\\
&\le C(|e_\phi|+|\rho_\phi|+|\nabla e_\phi|+|\nabla\rho_\phi|).
\end{aligned}
$$
It follows that
\begin{equation}
 |R_{\rm lift}|\le C\bigl(\nWi{e_\phi}{\Uha}+h^k\bigr)\|e_\phi\|_{H^1(\Gh)}\|e_u\|_{L^2(\Gh)} .
 \end{equation}

Then \eqref{eq:lift-remainder} follows from an inverse trace estimate. On the smooth exact surface, the tangential derivative maps $H^{1/2}$ to $H^{-1/2}$, and multiplication by the smooth coefficient $|\nabla\phi|^{-1}\nabla_\Gamma u$ is bounded on $H^{1/2}$. Thus $H^{-1/2}$--$H^{1/2}$ duality gives 
$$
 \bigl|\bigl(|\nabla\phi|^{-1}\nabla_{\Gamma[\phi]} u\cdot\nabla_{\Gamma[\phi]} e_\phi^\ell,e_u^\ell\bigr)_{\Gamma[\phi]}\bigr|
 \le C\|e_\phi^\ell\|_{H^{1/2}(\Gamma[\phi])}\|e_u^\ell\|_{H^{1/2}(\Gamma[\phi])}
$$
and the uniform trace and pullback bounds implied by the Lipschitz normal-graph property in \cref{subsec:continuation} complete the proof, for example,
$$
\|e_\phi^\ell\|_{H^{1/2}(\Gamma[\phi])}\leq C\|e_\phi\|_{H^{1/2}(\Gamma[\phi_h])}\leq C\nH{e_\phi}{\Uha} .
$$
\end{proof}

\Needspace{6\baselineskip}
\begin{proposition}[Bulk error controlled by the geometry error]
\label{prop:eu-final}
The following result holds for sufficiently small $h$ and $t\in[0,t^\ast)${\rm:}
\begin{equation}
 \nH{e_u}{\Oha}^2
 \le C\nH{e_\phi}{\Uha}^2+Ch^{2k} .                                 \label{eq:eu-final}
\end{equation}
\end{proposition}

\begin{proof}
Since $\Omega[\phi_h]$ is contained in $\Omega_h[\phi_h]$, the first term on the right-hand side of \eqref{eq:error-u-two-geometry} is bounded using \eqref{eq:r-shape-difference}:
$$
|(r,e_u)_{\Oh}|\leq Ch^k\|e_u\|_{H^1(\Omega_h[\phi_h])} . 
$$

The second term, $-(n(\nabla\phi_h)\cdot\nabla u,e_u)_{\Gh} $, is bounded using the decomposition
$$
 n(\nabla\phi_h)\cdot\nabla u
 =n(\nabla\phi) \cdot\nabla u +\bigl(n(\nabla\phi_h)-n(\nabla\phi)\bigr)\cdot\nabla u 
$$
and \cref{lem:bulk-terms}(a)--(c) together with \eqref{eq:ephi-Winf-from-H1}, which imply that
\begin{align*}
 |(n(\nabla\phi_h)\cdot\nabla u,e_u)_{\Gh}|
 &\le C\bigl(\|e_\phi\|_{H^1(U_h[\phi_h])}+h^k\bigr)\|e_u\|_{H^1(\Omega_h[\phi_h])} .
\end{align*}

The other terms on the right-hand side of \eqref{eq:error-u-two-geometry} are bounded by
$$
Ch^k\bigl(\|e_u\|_{H^1(\Omega_h[\phi_h])}+g_1(e_u,e_u)_{\Uha}^{1/2}\bigr),
$$
by \eqref{eq:interp-ass} and the face-jump interpolation bound following \eqref{eq:ghost-smooth-zero}. Combining these bounds, using \cref{lem:ghost-extension}, and applying Young's inequality proves \eqref{eq:eu-final}.
\end{proof}

\subsection{The reconstruction error}
\label{subsec:ephi-estimate}

\begin{proposition}[Reconstruction error]
\label{prop:ephi-final}
The following result holds for sufficiently small $h$ and $t\in[0,t^\ast)${\rm:}
\begin{equation}
 \nH{e_\phi}{U_h[\phi_h]}^2 
 \le C\nL{e_p}{\Uh}^2+C\nL{e_\eta}{\Uh}^2+Ch^{2k}.
 \label{eq:ephi-final}
\end{equation}
\end{proposition}

\begin{proof}
Taking $v_h=e_\phi$ in \eqref{eq:ephi-current-band} and using \eqref{eq:ephi-interp-residual} and Young's inequality yields \eqref{eq:ephi-final}, since $\nH{e_\phi}{U_h[\phi_h]}^2\le C_{\rm ext}\bigl(\nH{e_\phi}{\Uh}^2+g_1(e_\phi,e_\phi)_{\Uha}\bigr)$ by \cref{lem:ghost-extension}.  
\end{proof}

\subsection{Bounds on time derivatives of transport errors}
\label{subsec:material}

\begin{lemma}[Time derivatives of the transport errors]
\label{lem:dt-transport}
The following result holds for sufficiently small $h$ and $t\in[0,t^\ast)${\rm:} 
\begin{align*}
 &\nL{\partial_te_p}{\Uha}
 +\nL{\partial_te_\eta}{\Uha} \\
 &\quad\le C(h^{-1}\nL{e_p}{\Uha}+h^{-1}\nL{e_\eta}{\Uha}
 +\nH{e_u}{\Oha}+h^k)\\
 &\quad \le Ch^{k-1.1}.
\end{align*}
\end{lemma}

\begin{proof}
We take $v_h=\partial_te_p\in V_h|_{\Uha}$ in \eqref{eq:error-p-current-band} and move the transport term and the ghost term of $e_p$ to the right-hand side:
\begin{align*}
 &\nL{\partial_te_p}{\Uh}^2+g_0(\partial_te_p,\partial_te_p)_{\Uha}\\
 &\quad=-(\beta(u_h)n(\eta_h)\cdot\nabla e_p,\partial_te_p)_{\Uh}
 -(\partial_t\rho_p+\beta(u_h)n(\eta_h)\cdot\nabla\rho_p,\partial_te_p)_{\Uh}\\
 &\qquad-\bigl((\beta(u_h)n(\eta_h)-\beta(u)n(\eta))\cdot\nabla p,\partial_te_p\bigr)_{\Uh}\\
 &\qquad-g_0(e_p,\partial_te_p)_{\Uha}-g_0(\partial_t\rho_p+\rho_p,\partial_te_p)_{\Uha}\\
 &\quad\leq Ch^{-1}\nL{e_p}{\Uha}\nL{\partial_te_p}{\Uha}
 + Ch^k\nL{\partial_te_p}{\Uha} \\
 &\qquad + C(\nH{e_u}{\Oha}+\nL{e_\eta}{\Uha}+\nL{e_p}{\Uha})\nL{\partial_te_p}{\Uha} ,
\end{align*}
where the last inequality uses the inverse estimate in \eqref{eq:inverse-inequalities}. By \cref{lem:ghost-extension}, the left-hand side above is bounded below by $c\|\partial_te_p\|_{L^2(\Uha)}^2$ with $c>0$ independent of $h$. Hence
\begin{align*}
 &\nL{\partial_te_p}{\Uha} \le C(h^{-1}\nL{e_p}{\Uha}+\nL{e_\eta}{\Uha}
 +\nH{e_u}{\Oha}+h^k).
\end{align*}
For the vector equation, the additional source satisfies
$$
 \|\beta'(u_h)\nabla u_h|\eta_h|-\beta'(u)\nabla u|\eta|\|_{L^2(\Uh)}
 \le C\bigl(\nH{e_u}{\Oha}+\nL{e_\eta}{\Uha}+h^k\bigr).
$$
Indeed, separate the factor $\nabla(u_h-u)$ from the differences in $\beta'$ and $|\eta_h|$; the remaining coefficients are bounded, and interpolation controls $\nabla\rho_u$. Testing \eqref{eq:eeta-current-band} with $\zeta_h=\partial_te_\eta$ gives the same estimate for $\partial_te_\eta$. Together with \eqref{eq:continuation-conditions}, this proves the lemma.  
\end{proof}

\subsection{The geometry defect}
\label{subsec:geometry-defect}

We next estimate the defect caused by replacing the numerical band by the exact band. 

\begin{lemma}[Geometry defect]\label{lem:geometry-defect}
Let $w\in\{e_p,e_{\eta}\}$. For sufficiently small $h$ and $t\in[0,t^*)$,
\begin{align}
&\big|(\partial_tw+\beta(u)n(\nabla\phi)\cdot\nabla w,w)_{U[\phi_h]}-(\partial_tw+\beta(u)n(\nabla\phi)\cdot\nabla w,w)_{U[\phi]}\big|\nonumber\\
&\quad\le Ch^{-\frac{n+3}{4}}\big(\|e_\phi\|_{H^1(U_h)}+h^{k+1}\big)\big(\|\partial_tw\|_{L^2(U_h)}+\|\nabla w\|_{L^2(U_h)}\big)\|w\|_{L^2(U_h)}\label{eq:lem7.5_1}\\
&\quad\le C\big(\|e_p\|^2_{L^2(U_h)}+\|e_\eta\|^2_{L^2(U_h)}+\|e_u\|^2_{H^1(\Omega_h[\phi_h])}+h^{2k}\big),\label{eq:lem7.5_2}
\end{align}
where $U_h:=U_h[\phi_h]$. The second inequality requires $k\ge\frac{n+3}4+1.1$, so for $n\in\{2,3\}$ it holds for $k\ge3$.
\end{lemma}

\begin{proof}
Let $\chi^\theta:=\phi+\theta(\phi_h-\phi)$, $\theta\in[0,1]$, so that $U[\chi^0]=U[\phi]$, $U[\chi^1]=U[\phi_h]$ and $\partial_\theta\chi^\theta=\phi_h-\phi=\rho_\phi+e_\phi$. By \eqref{eq:interp-ass} and \eqref{eq:ephi-Winf-from-H1}, 
$$
|\nabla\chi^\theta-\nabla\phi|\le|\nabla\rho_\phi|+|\nabla e_\phi|\le Ch^k+h ,
$$ 
so $|\nabla\chi^\theta|\ge c/2$ on $U_h[\phi_h]$. By \eqref{eq:smallness}, $|\chi^\theta-\phi_h|\le|\phi-\phi_h|\le Ch^{1.1}$. Hence:
\begin{itemize}
\item
Every $x\in U[\chi^\theta]$ lies within the $Ch^{1.1}$-neighborhood of $U[\phi_h]$. 

\item \eqref{eq:band-buffer} gives $U[\chi^\theta]\subset U_h[\phi_h]$ with $\dist(\Gamma_\pm[\chi^\theta],\partial U_h[\phi_h])\ge ch$. 

\item The surfaces $S=\Gamma_\pm[\chi^\theta]$ are uniformly Lipschitz graphs over $\Gamma_\pm[\phi]$, within distance $Ch^{1.1}$. Each element $K$ meets $S$ in a set with $|S\cap K|\le Ch^{n-1}$.
\end{itemize}

\emph{Step 1: proof of \eqref{eq:lem7.5_1}.} 
For $w\in\{e_p,e_{\eta}\}$, by applying \eqref{eq:FTC-band} with the $\theta$-independent integrand $F:=(\partial_tw+\beta(u)n(\nabla\phi)\cdot\nabla w)\cdot w$ (which is elementwise smooth on $U_h[\phi_h]$), we have 
$$
(F,1)_{U[\phi_h]}-(F,1)_{U[\phi]}=\int_0^1\Big[\int_{\Gamma_-[\chi^\theta]}-\int_{\Gamma_+[\chi^\theta]}\Big]F\,\frac{\rho_\phi+e_\phi}{|\nabla\chi^\theta|}\dd s\dd\theta .
$$
By Hölder's inequality with exponents $(2,4,4)$ and $|\nabla\chi^\theta|\ge c/2$,
\begin{align}\label{F1-diff}
&\big|(F,1)_{U[\phi_h]}-(F,1)_{U[\phi]}\big| \notag\\
&\le C\sup_{\theta\in[0,1]}\sum_{S=\Gamma_\pm[\chi^\theta]}\|\partial_tw+\beta(u)n(\nabla\phi)\cdot\nabla w\|_{L^2(S)}\,\|w\|_{L^4(S)}\,\|\rho_\phi+e_\phi\|_{L^4(S)} .
\end{align}

For a piecewise polynomial $v$ on the elements $K\subset U_h[\phi_h]$ meeting $S$, the bounds $|S\cap K|\le Ch^{n-1}$ and $\|v\|_{L^\infty(K)}\le Ch^{-n/2}\|v\|_{L^2(K)}$ give
$$
\begin{aligned}
\|v\|^2_{L^2(S\cap K)}&\le|S\cap K|\,\|v\|^2_{L^\infty(K)}\le Ch^{-1}\|v\|^2_{L^2(K)},\\
\|v\|^4_{L^4(S\cap K)}&\le|S\cap K|\,\|v\|^4_{L^\infty(K)}\le Ch^{-(n+1)}\|v\|^4_{L^2(K)} .
\end{aligned}
$$
Summing over $K$, and using
$\sum_K\|v\|_{L^2(K)}^4\le(\sum_K\|v\|_{L^2(K)}^2)^2$
for the second estimate,
$$
\|v\|_{L^2(S)}\le Ch^{-1/2}\|v\|_{L^2(U_h[\phi_h])},\qquad \|v\|_{L^4(S)}\le Ch^{-\frac{n+1}4}\|v\|_{L^2(U_h[\phi_h])}.
$$
On each $I_m$ the active mesh is fixed, so $\partial_tw$ and $\nabla w$ are piecewise polynomial. Hence
\begin{align*}
\|\partial_tw+\beta(u)n(\nabla\phi)\cdot\nabla w\|_{L^2(S)}
&\le\|\partial_tw\|_{L^2(S)}+\|\beta(u)\|_{L^\infty(D)}\|\nabla w\|_{L^2(S)} \\
&\le Ch^{-1/2}\big(\|\partial_tw\|_{L^2(U_h[\phi_h])}+\|\nabla w\|_{L^2(U_h[\phi_h])}\big),
\end{align*}
and $\|w\|_{L^4(S)}\le Ch^{-(n+1)/4}\|w\|_{L^2(U_h[\phi_h])}$.

By \eqref{eq:interp-ass}, $\|\rho_\phi\|_{L^4(S)}\le |S|^{1/4}\|\rho_\phi\|_{L^\infty(S)}\le Ch^{k+1}$. For $e_\phi$, $S$ is a boundary component of $U[\chi^\theta]\subset U_h[\phi_h]$, a domain of width $\simeq\delta$ with uniformly Lipschitz boundary. The trace inequality $H^1(U[\chi^\theta])\to H^{1/2}(S)$ and the embedding $H^{1/2}(S)\hookrightarrow L^4(S)$ (valid in two and three dimensions) give the following estimate uniformly in $\theta$, $h$ and $t$:
$$
\|e_\phi\|_{L^4(S)}\le C\|e_\phi\|_{H^1(U_h[\phi_h])}.
$$
Substituting the above estimates into \eqref{F1-diff} proves \eqref{eq:lem7.5_1}, with $h^{-\frac{1}{2}}h^{-\frac{n+1}{4}}=h^{-\frac{n+3}{4}}$.

\emph{Step 2: proof of \eqref{eq:lem7.5_2}.} On $[0,t^*)$ we have $\|e_\phi\|_{H^1(U_h[\phi_h])}\le h^{k-0.1}$ by \eqref{eq:continuation-conditions}. The inverse inequality \eqref{eq:inverse-inequalities} gives $\|\nabla w\|_{L^2(U_h[\phi_h])}\le Ch^{-1}\|w\|_{L^2(U_h[\phi_h])}$. By \cref{lem:dt-transport}, 
$$
\begin{aligned}
\|\partial_tw\|_{L^2(U_h[\phi_h])}
 &\le C\bigl(h^{-1}\nL{e_p}{\Uha}+h^{-1}\nL{e_\eta}{\Uha}\bigr)\\
 &\quad+C\bigl(\|e_u\|_{H^1(\Omega_h[\phi_h])}+h^k\bigr).
\end{aligned}
$$
With $\nL{w}{\Uha} \le \nL{e_p}{\Uha} + \nL{e_\eta}{\Uha}$, \eqref{eq:lem7.5_1} gives
$$
\begin{aligned}
&\big|(\partial_tw+\beta(u)n(\nabla\phi)\cdot\nabla w,w)_{U[\phi_h]}-(\partial_tw+\beta(u)n(\nabla\phi)\cdot\nabla w,w)_{U[\phi]}\big|\\
&\le Ch^{k-1.1-\frac{n+3}4}\big(\nL{e_p}{\Uha}+\nL{e_\eta}{\Uha}\big)^2\\
&\quad+Ch^{k-0.1-\frac{n+3}4}\big(\|e_u\|_{H^1(\Omega_h[\phi_h])}+h^k\big)\big(\nL{e_p}{\Uha}+\nL{e_\eta}{\Uha}\big).
\end{aligned}
$$
Both powers of $h$ are nonnegative when $k\ge1.1+\frac{n+3}4$, which holds for $k\ge3$ and $n\le3$. Young's inequality then gives \eqref{eq:lem7.5_2}.
\end{proof}

\subsection{The transport errors on the exact band}
\label{subsec:transport-estimate}

\begin{lemma}[Reynolds identity on the exact band]
\label{lem:reynolds}
On an interval $I_m$, let $v\in C^1(I_m;V_h|_{\Uha})$ and suppose that the exact band is compactly contained in $\Uha$. Then, away from the interval endpoints,
\begin{equation}
 (\partial_tv+\beta(u)n(\nabla\phi)\cdot\nabla v,v)_{\Ub{\phi}}
 =\frac12\frac{\dd}{\dd t}\nL{v}{\Ub{\phi}}^2
 -\frac12\bigl(\nabla\cdot(\beta(u)n(\nabla\phi)),v^2\bigr)_{\Ub{\phi}} .
 \label{eq:reynolds-identity}
\end{equation}
\end{lemma}

\begin{proof}
By \eqref{eq:HJ}, $\partial_t\phi=-\beta(u)n(\nabla\phi)\cdot\nabla\phi$ in $\Nd(t)$, so that \eqref{eq:reynolds-band} reads
$$
 \frac{\dd}{\dd t}\nL{v}{\Ub{\phi}}^2
 =2(\partial_tv,v)_{\Ub{\phi}}
 +\int_{\Gamma_+[\phi]}v^2\,\beta(u)\,\dd s
 -\int_{\Gamma_-[\phi]}v^2\,\beta(u)\,\dd s ,
$$
and the outward unit normal $\nu$ of $\Ub{\phi}$ is $n(\nabla\phi)$ on $\Gamma_+[\phi]$ and $-n(\nabla\phi)$ on $\Gamma_-[\phi]$, so the two surface integrals equal $\int_{\partial\Ub{\phi}}v^2\beta(u)n(\nabla\phi)\cdot\nu\dd s$. The divergence theorem, applied elementwise with cancellation across interior faces, gives
$$
 2(\beta(u)n(\nabla\phi)\cdot\nabla v,v)_{\Ub{\phi}}
 =\int_{\partial\Ub{\phi}}v^2\beta(u)n(\nabla\phi)\cdot\nu\,\dd s
 -(\nabla\cdot(\beta(u)n(\nabla\phi)),v^2)_{\Ub{\phi}}.
$$
Eliminating the boundary integral between the two identities gives \eqref{eq:reynolds-identity}.
\end{proof}

We measure the errors $e_p$ and $e_\eta$ in the following norms: 
\begin{align}
\begin{aligned}
 \|e_p\|_{E}^2&:= \frac12 \nL{e_p}{\Ub{\phi}}^2+\frac12 g_0(e_p,e_p)_{\Uha}, \\
 \|e_\eta\|_{E}^2&:=\frac12 \nL{e_\eta}{\Ub{\phi}}^2+\frac12 g_0(e_\eta,e_\eta)_{\Uha}. 
\end{aligned}
 \label{eq:energy-def}
\end{align}
By \cref{lem:ghost-extension}, applied with the band $\Ub{\phi}$ (the chain condition is verified in \cref{appendix:path}),
\begin{equation}
 \nL{e_p}{\Uha}\le C\|e_p\|_E,
 \qquad
 \nL{e_\eta}{\Uha}\le C\|e_\eta\|_E.
 \label{eq:energy-lower}
\end{equation}

\Needspace{6\baselineskip}
\begin{proposition}[Transport errors]
\label{prop:transport}
The following result holds for sufficiently small $h$ and every time interval $I_m\cap[0,t^\ast)${\rm:}
\begin{align}
 \frac{\dd}{\dd t} \|e_p\|_{E}^2 +g_0(e_p,e_p)_{\Uha}
 &\le C( \|e_p\|_{E}^2+ \|e_\eta\|_{E}^2+ \nH{e_u}{\Oha}^2 +h^{2k}),
 \label{eq:ep-differential}\\
 \frac{\dd}{\dd t} \|e_\eta\|_{E}^2 +g_0(e_\eta,e_\eta)_{\Uha}
 &\le C( \|e_p\|_{E}^2+ \|e_\eta\|_{E}^2 + \nH{e_u}{\Oha}^2 + h^{2k}).
 \label{eq:eeta-differential}
\end{align}
\end{proposition}

\begin{proof}
We take $v_h=e_p$ in \eqref{eq:error-p-current-band} and rewrite the first term on the left-hand side as
\begin{align*}
 &(\partial_te_p+\beta(u_h)n(\eta_h)\cdot\nabla e_p,e_p)_{\Uh}\\
 &\quad=(\partial_te_p+\beta(u)n(\eta)\cdot\nabla e_p,e_p)_{\Ub{\phi}} +\bigl((\beta(u_h)n(\eta_h)-\beta(u)n(\eta))\cdot\nabla e_p,e_p\bigr)_{\Uh} \\
 &\qquad + (\partial_te_p+\beta(u)n(\eta)\cdot\nabla e_p,e_p)_{\Uh}
 -(\partial_te_p +\beta(u)n(\eta)\cdot\nabla e_p,e_p)_{\Ub{\phi}} .
\end{align*}
\cref{lem:reynolds} with $v=e_p$ turns the first term on the right-hand side of the last equality into $\frac12\frac{\dd}{\dd t}\nL{e_p}{\Ub{\phi}}^2-\frac12(\nabla\cdot(\beta(u)n(\nabla\phi)),e_p^2)_{\Ub{\phi}}$. Since $g_0(\partial_te_p,e_p)_{\Uha} = \frac12 \frac{\dd}{\dd t}g_0(e_p,e_p)_{\Uha}$ on a time interval $I_m$, it follows that 
\begin{align}
&\frac{\dd}{\dd t} \|e_p\|_{E}^2 + g_0(e_p,e_p)_{\Uha}\notag\\
&=\frac12\bigl(\nabla\cdot(\beta(u)n(\eta)),e_p^2\bigr)_{\Ub{\phi}} \notag\\
  &\quad-\bigl((\beta(u_h)n(\eta_h)-\beta(u)n(\eta))\cdot\nabla e_p,e_p\bigr)_{\Uh} \notag\\
  &\quad - [(\partial_te_p+\beta(u)n(\eta)\cdot\nabla e_p,e_p)_{\Uh}
 -(\partial_te_p +\beta(u)n(\eta)\cdot\nabla e_p,e_p)_{\Ub{\phi}} ]\notag\\
 &\quad-(\partial_t\rho_p+\beta(u_h)n(\eta_h)\cdot\nabla\rho_p,e_p)_{\Uh} \notag\\
 &\quad-\bigl((\beta(u_h)n(\eta_h)-\beta(u)n(\eta))\cdot\nabla p,e_p\bigr)_{\Uh} \notag\\
 &\quad-g_0(\partial_t\rho_p+\rho_p,e_p)_{\Uha} \notag\\
 &=:J_1+J_2+J_3+J_4+J_5+J_6 ,
\label{eq:ep-energy-identity}
\end{align}
where 
$$
|J_1|+|J_5|\leq C(\|e_p\|_{L^2(U_h[\phi_h])} + \|e_u\|_{L^2(U_h[\phi_h])} + \|e_\eta\|_{L^2(U_h[\phi_h])} + h^k )\|e_p\|_{L^2(U_h[\phi_h])} ,
$$ 
and \eqref{eq:smallness} with the inverse estimate \eqref{eq:inverse-inequalities} implies that 
$$
|J_2|
\leq Ch^{-1}(\|e_u\|_{L^\infty(\Omega_h[\phi_h])} + \|e_\eta\|_{L^\infty(U_h[\phi_h])} + h^{k+1}) \|e_p\|_{L^2(U_h[\phi_h])}^2
\leq C\|e_p\|_{L^2(U_h[\phi_h])}^2 .
$$
The estimate of $J_3$ follows from \cref{lem:geometry-defect}: 
\begin{equation}
 |J_3|\le C\bigl(\nL{e_p}{\Uha}^2+\nL{e_\eta}{\Uha}^2 + \nH{e_u}{\Oha}^2+h^{2k}\bigr).
 \label{eq:J2}
\end{equation}
By \eqref{eq:Dm-rho-bound}, \eqref{eq:ghost-rho-bound} and the inverse trace estimate $g_0(e_p,e_p)_{\Uha}\le C\nL{e_p}{\Uha}^2$, we have 
\begin{equation}
 |J_4|+|J_6|\le Ch^k \nL{e_p}{\Uha} .
 \label{eq:J4-J6}
\end{equation}
Collecting these bounds yields \eqref{eq:ep-differential}. For \eqref{eq:eeta-differential}, apply the same argument componentwise
and use the source bound in the proof of \cref{lem:dt-transport}, paired with $e_\eta$, followed by Young's inequality.
\end{proof}

\subsection{Transfer between active meshes}
\label{subsec:transfer}

\begin{lemma}[Transfer]
\label{lem:transfer}
Let $t_{m+1}\le t^\ast$ be an update time {\rm(}by the definition of $t^\ast$, and by \cref{sec:closing} if $t_{m+1}=t^\ast$, \eqref{eq:continuation-conditions} holds for the left limit at $t_{m+1}$ on the old active mesh{\rm)}. Then
\begin{align}\label{eq:transfer-energy}
\begin{aligned}
\|e_p(t_{m+1}^+)\|_{E}^2 &\leq \bigl(1+Ch^{k-0.1-\frac{n+3}{4}}\bigr) \|e_p(t_{m+1}^-)\|_{E}^2+Ch^{2k+2} ,\\
\|e_\eta(t_{m+1}^+)\|_{E}^2 &\leq \bigl(1+Ch^{k-0.1-\frac{n+3}{4}}\bigr) \|e_\eta(t_{m+1}^-)\|_{E}^2+Ch^{2k+2} ,
\end{aligned}
\end{align}
where
$$
\|w(t_{m+1}^\pm)\|_{E}^2:=\tfrac12\nL{w(t_{m+1}^\pm)}{\Ub{\phi}}^2
+\tfrac12g_0\bigl(w(t_{m+1}^\pm),w(t_{m+1}^\pm)\bigr)_{U_h^\pm}
$$
for $w\in\{e_p,e_\eta\}$, and $U_h^-$ and $U_h^+$ are the active meshes used before and after the update, built from $\phi_h(t_m^-)$ and $\phi_h(t_{m+1}^-)$, respectively.
\end{lemma}

\begin{proof}
Write $p_h^\pm=p_h(t_{m+1}^\pm)$, $\phi_h^\pm:=\phi_h(t_{m+1}^\pm)$, $e_p^\pm=e_p(t_{m+1}^\pm)$, $e_\phi^-:=e_\phi(t_{m+1}^-)$, $U^-:=\Ub{\phi_h^-}$, and $U_h^-$, $U_h^+$ for the old and new active meshes, built from $\phi_h(t_m^-)$ and $\phi_h^-$. In particular, $U_h^+$ is constructed from $\phi_h^-$ before $\phi_h^+$ is reconstructed. The exact band and the bands $U[\chi^\theta]$ of the homotopy from $\phi$ to $\phi_h^-$ lie in both active meshes, and the chain condition of \cref{appendix:path} holds on both active meshes. These facts use only the pretransfer geometry. Then, by \eqref{eq:stabilized-extension}, 
$$
 (p_h^+,v_h)_{U^-}+g_0(p_h^+,v_h)_{U_h^+}=(p_h^-,v_h)_{U^-}
 \quad\forall v_h\in V_h|_{U_h^+} .
$$
Since $g_0(p,\cdot)_{U_h^+}=0$ by \eqref{eq:ghost-smooth-zero}, subtracting $(I_hp,v_h)_{U^-}+g_0(I_hp,v_h)_{U_h^+}$ from both
sides of the equation above gives
$$
(e_p^+,v_h)_{U^-}+g_0(e_p^+,v_h)_{U_h^+}=(e_p^-,v_h)_{U^-}-g_0(\rho_p,v_h)_{U_h^+} 
\quad\forall v_h\in V_h|_{U_h^+} .
$$
Taking $v_h=e_p^+$, using $g_0(\rho_p,\rho_p)_{U_h^+}\le Ch^{2k+2}$ by the ghost interpolation bound following \eqref{eq:ghost-smooth-zero}, and Young's inequality, we have 
\begin{equation}
 \nL{e_p^+}{U^-}^2+g_0(e_p^+,e_p^+)_{U_h^+}\le\nL{e_p^-}{U^-}^2+Ch^{2k+2}.
 \label{eq:transfer-step1}
\end{equation}

With $\chi^\theta:=\phi+\theta(\phi_h^--\phi)$, applying \eqref{eq:FTC-band} with the integrand $(e_p^\pm)^2$ gives
\begin{align}\label{epU-Uphi}
\Big| \nL{e_p^\pm}{U^-}^2-\nL{e_p^\pm}{\Ub{\phi}}^2 \Big| 
& = \Big| \int_0^1\Bigl[\int_{\Gamma_-[\chi^\theta]}-\int_{\Gamma_+[\chi^\theta]}\Bigr]
 (e_p^\pm)^2\frac{\rho_\phi+e_\phi^-}{|\nabla\chi^\theta|}\,\dd s\,\dd\theta \Big| \notag\\
&\leq Ch^{k-0.1-\frac{n+3}{4}}\nL{e_p^\pm}{U_h^\pm}^2 \notag\\
&\leq Ch^{k-0.1-\frac{n+3}{4}}\|e_p^\pm\|_{E}^2 ,
\end{align}
where the first inequality uses the $L^\infty$ bound on $e_\phi^-$ from \eqref{eq:discrete-sobolev} and \eqref{eq:continuation-conditions}, and the inverse trace estimate $\nL{e_p^\pm}{S}^2\leq Ch^{-1}\nL{e_p^\pm}{U_h^\pm}^2$ for $S=\Gamma_\pm[\chi^\theta]$; the last inequality follows from \eqref{eq:energy-lower}. 

Apply \eqref{epU-Uphi} separately on the two active meshes, bounding each active-mesh norm by its own stabilized energy. Together with \eqref{eq:transfer-step1}, and discarding the nonnegative old ghost term, this gives
$$
 \|e_p^+\|_E^2\le\|e_p^-\|_E^2
 +Ch^{k-0.1-(n+3)/4}
       \bigl(\|e_p^+\|_E^2+\|e_p^-\|_E^2\bigr)+Ch^{2k+2}.
$$
Since $k-0.1-(n+3)/4>0$, absorption proves \eqref{eq:transfer-energy}. Apply the same argument componentwise to $e_\eta$.
\end{proof}

\subsection{The coupled estimate}
\label{subsec:coupled}

\begin{proposition}[Coupled estimate]
\label{prop:coupled}
The following result holds for sufficiently small $h$ and $t\in[0,t^\ast)${\rm:}
\begin{equation}
 \nL{e_p}{\Uha}+\nL{e_\eta}{\Uha}+\nH{e_\phi}{\Uha}+\nH{e_u}{\Oha}
 \le C_0h^{k} . 
 \label{eq:coupled-pointwise}
\end{equation}
\end{proposition}

\begin{proof}
Adding \eqref{eq:ep-differential} and \eqref{eq:eeta-differential}, using \cref{prop:eu-final,prop:ephi-final} and \eqref{eq:energy-lower},
gives
$$
 \frac{\dd}{\dd t}\bigl(\|e_p\|_E^2+\|e_\eta\|_E^2\bigr)
 \le C\bigl(\|e_p\|_E^2+\|e_\eta\|_E^2\bigr)+Ch^{2k}
$$
inside every interval $I_m$. Hence, for $t_m<t<\min\{t_{m+1},t^\ast\}$, Gronwall's inequality yields
$$
 \|e_p(t)\|_E^2+\|e_\eta(t)\|_E^2
 \le e^{C(t-t_m)}\bigl(\|e_p(t_m^+)\|_E^2
       +\|e_\eta(t_m^+)\|_E^2+C(t-t_m)h^{2k}\bigr).
$$
At an update $t_{m+1}<t^\ast$, pass to the left limit $t_{m+1}^-$ on the old active meshes and use \cref{lem:transfer}. Since $k-0.1-(n+3)/4>1$ and $h\le C(t_{m+1}-t_m)$ by \eqref{time-length-min},
\begin{align*}
 &\|e_p(t_{m+1}^+)\|_E^2+\|e_\eta(t_{m+1}^+)\|_E^2\\
 &\quad\le(1+Ch)e^{C(t_{m+1}-t_m)}
 \bigl(\|e_p(t_m^+)\|_E^2+\|e_\eta(t_m^+)\|_E^2
       +C(t_{m+1}-t_m)h^{2k}\bigr)+Ch^{2k+2}\\
 &\quad\le e^{C'(t_{m+1}-t_m)}
 \bigl(\|e_p(t_m^+)\|_E^2+\|e_\eta(t_m^+)\|_E^2
       +C(t_{m+1}-t_m)h^{2k}\bigr).
\end{align*}
For any $t<t^\ast$, iterate this estimate over the preceding updates and apply the first Gronwall bound on the final partial interval.
The initial transport errors vanish and the elapsed times sum to at most $T$, so
$$
 \sup_{0\le t<t^\ast}\bigl(\|e_p(t)\|_E^2+\|e_\eta(t)\|_E^2\bigr)
 \le Ch^{2k}.
$$
This includes both one-sided values at each update. Then \eqref{eq:energy-lower} and \cref{prop:eu-final,prop:ephi-final} prove \eqref{eq:coupled-pointwise}.
\end{proof}

\subsection{Completion of the proof}

Under the weaker bounds \eqref{eq:continuation-conditions}, we have established the improved error estimate of \cref{prop:coupled}
for $t\in[0,t^\ast)$. This improvement allows us to close the continuation argument and conclude that $t^\ast=T$, which completes
the proof of \cref{thm:main}. This argument is standard in the error analysis of semidiscrete FEMs for nonlinear PDEs \cite{FernandezLiOlshanskii2026,KLLP2017,KLL2019,KLL2021}. The detailed argument for the present scheme is given in \cref{sec:closure}.

\section{Numerical results}
\label{sec:num}

In this section, we test the convergence rates of the proposed method on an expanding ellipse, and present numerical simulation of a PDE-driven translating, expanding and rotating flower-shape domain. 

We use our own Python code (NumPy, SciPy) on meshes of squares of side $h$, each cut into two triangles. On a cut triangle we integrate over $\{\phi_h<c\}$. For the ellipse we replace the curve $\{\phi_h=c\}$ by a polyline with vertices on the curve {and use the signed fan triangulation of \cite[Section~4.2]{BCDE2021}, with Duffy quadrature on each triangle}. For the flower we use the quadrature of Saye \cite{Saye2015}{, implemented with Algoim (order $10$) on three quadrilateral subcells per triangle}. 

\subsection{Convergence test on an expanding ellipse}
\label{subsec:num-ellipse}

Let $E$ be the filled ellipse with semi-axes $0.6$ and $0.25$, rotated by the angle $0.3$ and centered at $c=(0.0137,-0.0219)$. The
boundary $\Gamma(t)$ is the parallel curve of $E$ at distance $R(t)=0.7+0.3t$, and $\phi(x,t)=\dist(x,E)-R(t)$.

We take $\beta(s)=s$ and $u=0.3+0.5\,\chi\, g$, where $g(x)=\cos(2(x_1-c_1))\cos(2(x_2-c_2))$ and $\chi=(1-S)_+^6$. Here $S$ is
the quadratic form of the ellipse with semi-axes $0.68$ and $0.33$ and the same center and rotation as $E$. The bump
$\chi g$ vanishes on $\Nd(t)$. Hence $u=0.3$ near $\Gamma(t)$, the boundary moves with normal speed $0.3$, and $f=u-\Delta u$.

The domain is $D=(-2,2)^2$ and $T=0.5$. We use uniform meshes with $h=4/n$, $n\in\{56,64,80,96,112,128\}$, and take $k=1,2,3$, $\kappa=2$, $\delta=0.3$, $c_{\rm buf}=4$, $\gamma_{0,\ell}=1/(\ell!)^2$, $\gamma_{1,\ell}=0.05/(\ell!)^2$ and $\tau=h/2$.
Since $p$ is affine in $t$ and $\eta$ does not depend on $t$, the time stepping adds no consistency error.

\begin{figure}[!htbp]
\centering
\begin{minipage}[t]{0.30\textwidth}
\centering
\includegraphics[width=\linewidth]{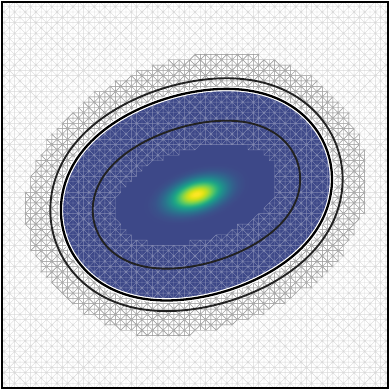}\\[2pt]
{\small (a) $t=0$}
\end{minipage}\hspace{0.06\textwidth}
\begin{minipage}[t]{0.30\textwidth}
\centering
\includegraphics[width=\linewidth]{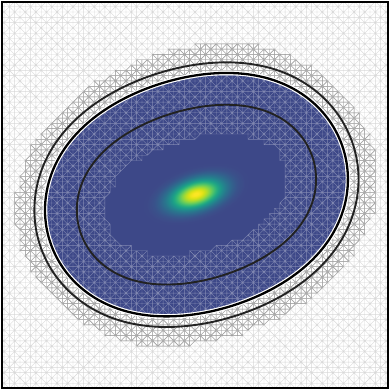}\\[2pt]
{\small (b) $t=0.5$}
\end{minipage}

\vspace{6pt}
\includegraphics[width=0.45\textwidth]{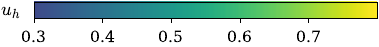}
\caption{$u_h$ on $\Omega[\phi_h]$ computed with $k=3$ and $h=0.05$. Black line:
$\Gamma[\phi_h]$; dashed lines: $\partial U[\phi_h]$; gray: elements of $U_h[\phi_h]$.}
\label{fig:ellipse-evolution}
\end{figure}

\begin{figure}[!htbp]
\centering
\includegraphics[width=0.40\textwidth]{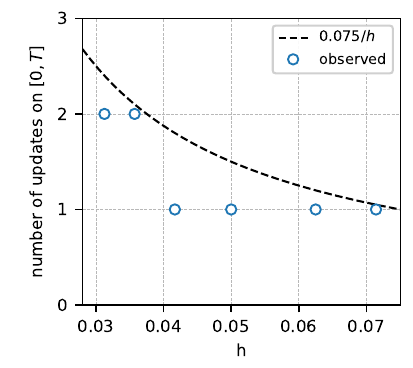}
\caption{Number of active-mesh updates on $[0,T]$ for the ellipse.}
\label{fig:ellipse-updates}
\end{figure}

\begin{figure}[!htbp]
\centering
\begin{minipage}[t]{0.42\textwidth}
\centering
\includegraphics[width=\linewidth]{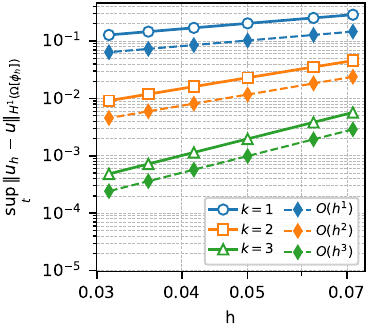}\\[1pt]
{\small (a) $u_h-u$ on $\Omega[\phi_h]$, $k=1,2,3$}
\end{minipage}\hfill
\begin{minipage}[t]{0.42\textwidth}
\centering
\includegraphics[width=0.96\linewidth]{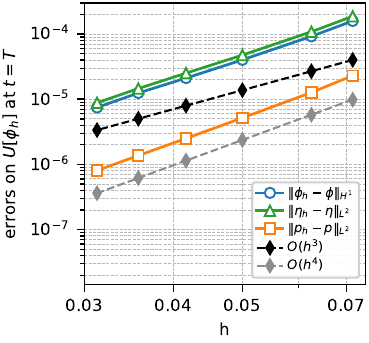}\\[1pt]
{\small (b) errors for $\phi_h$, $p_h$, $\eta_h$ on $U[\phi_h]$, $k=3$}
\end{minipage}
\vspace{-2pt}
\caption{Errors for the ellipse.}
\label{fig:ellipse-u}
\end{figure}

\Cref{fig:ellipse-evolution} shows $u_h$, the band and the active mesh
for $k=3$ and $h=0.05$ at $t=0$ and $t=T$. The number of updates depends on $h$. Up to $T=0.5$ the active mesh is updated once for $n\le96$ and twice
for $n=112,128$; see \cref{fig:ellipse-updates}. The band
moves by about $\tfrac12c_{\rm buf}h$ between two updates, so the number of
updates grows like $2TV/(c_{\rm buf}h)=0.075/h$, where $V=0.3$ is the normal speed.

\Cref{fig:ellipse-u}(a) shows the error of $u_h$ on the active mesh, measured as the maximum over all time steps. \Cref{fig:ellipse-u}(b) shows the errors of $\phi_h$, $p_h$, and $\eta_h$ on the band $\Uh$ at $t=T$ for $k=3$. The numerical results indicate $k$th-order convergence for $k=1,2,3$, although our error analysis establishes this result only for $k\geq 3$. This suggests that the restriction $k\geq 3$ is a technical artifact of the present error-analysis framework rather than an intrinsic limitation of the method.

\subsection{A moving five-petal flower}
\label{subsec:num-petals}

Let $C$ be the curve $r=0.8\,(1+0.07\cos(5\theta+0.4))$ in polar coordinates around $c=(0.0137,-0.0219)$. The boundary $\Gamma(t)$ is the parallel curve of $C$ at distance $0.32+0.05t$, translated by $(0.1,0.04)\,t$ and rotated by the angle $0.3t$ about the moving center.

We take $\beta(s)=s$. Near $\Gamma(t)$, $u$ equals the normal velocity of $\Gamma(t)$; inside, $u$ has a smooth bump. The source is $f=u-\Delta u$. The domain is $D=(-1.8,1.8)^2$, $h=3.6/n$, $T=0.5$, $\kappa=2$, $\delta=0.149$ and $c_{\rm buf}=4$. 

\begin{figure}[h]
\centering
\centering
\begin{minipage}[t]{0.27\textwidth}
\centering
\includegraphics[width=\linewidth]{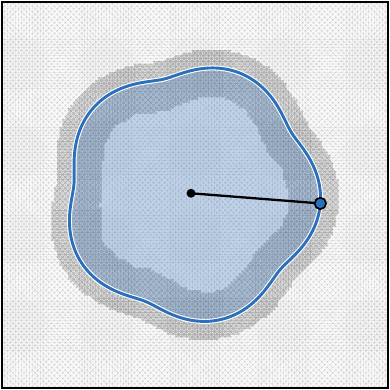}\\[2pt]
{\small (a) $t=0$}
\end{minipage}\hfill
\begin{minipage}[t]{0.27\textwidth}
\centering
\includegraphics[width=\linewidth]{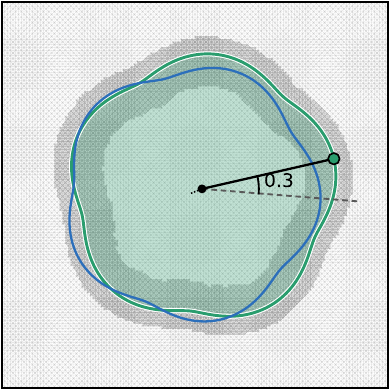}\\[2pt]
{\small (b) $t=1$}
\end{minipage}\hfill
\begin{minipage}[t]{0.27\textwidth}
\centering
\includegraphics[width=\linewidth]{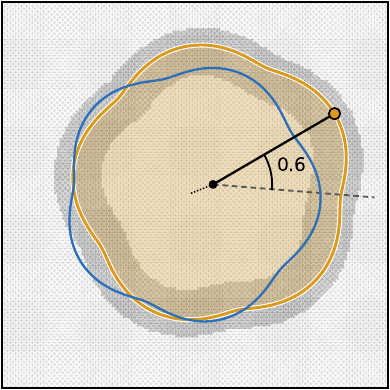}\\[2pt]
{\small (c) $t=2$}
\end{minipage}
\caption{Five-petal flower, $k=3$, $h=0.025$, on the fixed background mesh. Colored:
$\Omega[\phi_h]$; gray: $U_h$. Black: center to a petal tip; dashed: its direction at
$t=0$, with the angle turned. In (b), (c), blue: $\Gamma[\phi_h]$ at $t=0$; dotted: path of the center.}
\label{fig:petals-curves}
\end{figure}

We present numerical simulation up to $T=2$ with $k=3$, $h=0.025$, $\tau=0.0125$ and $\delta=0.145$. \Cref{fig:petals-curves} shows the domain and the active mesh at $t=0,1,2$ on the same background mesh: the flower moves to the right, turns and grows, and the active mesh follows it. The active mesh is rebuilt eight times in the finest mesh. The proposed method successfully captures the motion of the translating, expanding and rotating flower-shape domain.

\section{Conclusion}
\label{sec:conclusion}

We have proposed and analyzed an unfitted CutFEM for the quasistatic problem \eqref{eq:intro-model}, in which the boundary moves with normal velocity $\beta(u)$ depending on the solution $u$ of an elliptic PDE on the moving domain. The domain is given by a level-set function $\phi$. We evolve $p=\phi$ and $\eta=\nabla\phi$ by transport equations on a one-sided band and reconstruct $\phi$ from $\phi-\Delta\phi=p-\nabla\cdot\eta$. The active meshes are rebuilt only when the band comes close to their boundary, and the transport variables are then transferred by stabilized extension. For smooth evolution and $k\ge3$, we have proved $k$th-order convergence in $H^1$ for $u$ and $\phi$ and in $L^2$ for $p$ and $\eta$. In the numerical experiments we have observed $k$th-order convergence for all $k=1,2,3$. 

\section*{Acknowledgement}
\label{sec:acknowledgement}
The authors would like to thank Prof. Rong Tang for implementing the numerical scheme and providing the numerical results presented in Section~8. The corresponding code will be made publicly available on GitHub. The authors used Claude Fable 5.1 to assist with text editing, and take full responsibility for the accuracy, originality, and integrity of the manuscript.

\newpage
\appendix

\section{Full discretization}
\label{subsec:fd-scheme}
In computations we use a semi-implicit backward Euler scheme. Let
$$
 0=t_0^{\rm fd}<t_1^{\rm fd}<\cdots<t_J^{\rm fd}=T,
 \qquad
 \tau_j:=t_{j+1}^{\rm fd}-t_j^{\rm fd}\le\tau,
 \quad j=0,\ldots,J-1,
$$
where $\tau>0$ is a prescribed maximal step size. We denote the numerical band and bulk domain at time $t_j^{\rm fd}$ by 
$$ U^j:=U[\phi_h^j] \quad\mbox{and}\quad \Omega^j:=\Omega[\phi_h^j], $$
and their active meshes by $U_h^j$ and $\Omega_h^j$, with $U_h^j\subset\Omega_h^j$. These meshes remain fixed between mesh updates. 

The initial values $p_h^0,\eta_h^0,\phi_h^0$ and the initial active meshes are those of \cref{subsec:scheme}. Given $(p_h^{j-1},\eta_h^{j-1},\phi_h^{j-1},u_h^{j-1})$, we first compute provisional transport and reconstruction variables
$$
 p_h^{j,-},\phi_h^{j,-}\in V_h|_{U_h^{j-1}},
 \qquad
 \eta_h^{j,-}\in\bigl(V_h|_{U_h^{j-1}}\bigr)^n.
$$
The superscript $j,-$ indicates that these quantities are computed on the current active mesh, before any mesh update at $t_j^{\rm fd}$, as $w(t_m^-)$ in \cref{subsec:active}. They satisfy
\begin{align}
 &\Bigl(
   \frac{p_h^{j,-}-p_h^{j-1}}{\tau_{j-1}}
   +\beta(u_h^{j-1})n(\eta_h^{j-1})\cdot\nabla p_h^{j,-},
   v_h
 \Bigr)_{U^{j-1}}
 \nonumber\\
 &\qquad
 +g_0\Bigl(
   \frac{p_h^{j,-}-p_h^{j-1}}{\tau_{j-1}}+p_h^{j,-},
   v_h
 \Bigr)_{U_h^{j-1}}
 =0,
 \label{eq:fd-p}\\
 &\Bigl(
   \frac{\eta_h^{j,-}-\eta_h^{j-1}}{\tau_{j-1}}
   +\beta(u_h^{j-1})n(\eta_h^{j-1})\cdot\nabla\eta_h^{j,-},
   \zeta_h
 \Bigr)_{U^{j-1}}
 \nonumber\\
 &\qquad
 +g_0\Bigl(
   \frac{\eta_h^{j,-}-\eta_h^{j-1}}{\tau_{j-1}}+\eta_h^{j,-},
   \zeta_h
 \Bigr)_{U_h^{j-1}}
 \nonumber\\
 &\qquad
 =-\bigl(
   \beta'(u_h^{j-1})\nabla u_h^{j-1}|\eta_h^{j-1}|,
   \zeta_h
 \bigr)_{U^{j-1}},
 \label{eq:fd-eta}\\
 &(\phi_h^{j,-},v_h)_{U^{j-1}}
 +(\nabla\phi_h^{j,-},\nabla v_h)_{U^{j-1}}
 +g_1(\phi_h^{j,-},v_h)_{U_h^{j-1}}
 \nonumber\\
 &\qquad
 =(p_h^{j,-},v_h)_{U^{j-1}}
 +(\eta_h^{j,-},\nabla v_h)_{U^{j-1}},
 \label{eq:fd-rec}
\end{align}
for all $v_h\in V_h|_{U_h^{j-1}}$ and $\zeta_h\in\bigl(V_h|_{U_h^{j-1}}\bigr)^n$. All integrals are over the known band $U^{j-1}$, so the transport steps and the reconstruction are linear. We require that $|\eta_h^{j-1}|$ is bounded away from zero on $U_h^{j-1}$ and that the elements of $U_h^{j-1}$ and the band $U^{j-1}$ satisfy the chain condition, so that the stabilized mass form below is positive definite by \cref{lem:ghost-extension}, and choose $\tau_{j-1}$ so that
\begin{align*}
 &\tau_{j-1}\bigl|
 (\beta(u_h^{j-1})n(\eta_h^{j-1})\cdot\nabla v_h,v_h)_{U^{j-1}}
 \bigr| \le\tfrac12\bigl(
 \|v_h\|_{L^2(U^{j-1})}^2+g_0(v_h,v_h)_{U_h^{j-1}}\bigr) 
\end{align*}
for all $v_h\in V_h|_{U_h^{j-1}}$. 
This condition uses only the current numerical state and can be checked from the assembled forms. The transport bilinear form is then bounded below by $(2\tau_{j-1})^{-1}$ times the stabilized mass form, so both transport systems are invertible. The bilinear form on the left-hand side of \eqref{eq:fd-rec} is symmetric and, by the second estimate in \eqref{eq:proved-ghost-extension}, positive definite, so \eqref{eq:fd-rec} has a unique solution. 

After computing $\phi_h^{j,-}$, set $U^{j,-}:=U[\phi_h^{j,-}]$. Before any transfer, check that the provisional level surfaces are
nondegenerate and have the inherited component and side classification, that $\inf_{x\in U_h^{j-1}}|\eta_h^{j,-}(x)|>0$, and that
\begin{equation}
 U^{j,-}\Subset U_h^{j-1}.
 \label{eq:fd-update-test0}
\end{equation}
A failed trial is discarded and recomputed from the preceding accepted state with a smaller $\tau_{j-1}$; its data are not transferred. For an admissible trial, check the analogue of \eqref{eq:band-buffer},
\begin{align}\label{eq:fd-update-test}
\begin{aligned}
 \bigl\{x\in D:
        \dist(x,U^{j,-})\le\tfrac12c_{\rm buf}h
 \bigr\}
 &\Subset U_h^{j-1}, \\
 \bigl\{x\in D:
        \dist(x,U^{j,-})\le 2c_{\rm buf}h
 \bigr\} 
 &\supset U_h^{j-1}.
\end{aligned}
\end{align}
If \eqref{eq:fd-update-test} holds, we keep $ U_h^j:=U_h^{j-1}$ and $ \Omega_h^j:=\Omega_h^{j-1}$, and set
$
 (p_h^j,\eta_h^j,\phi_h^j)
 :=
 (p_h^{j,-},\eta_h^{j,-},\phi_h^{j,-}).
$

If \eqref{eq:fd-update-test} fails, we rebuild the active meshes from $\phi_h^{j,-}$ by \eqref{eq:buffer-elements}--\eqref{eq:numerical-bulk-support}. The provisional transport variables are then transferred to the new active mesh $U_h^j$ by the stabilized extensions,
$$
 p_h^j:=E_{U_h^j}p_h^{j,-},
 \qquad
 \eta_h^j:=E_{U_h^j}\eta_h^{j,-}.
$$
Both extensions use the band $U^{j,-}$, as in \eqref{eq:transfer}. Using these transferred variables, we reconstruct $\phi_h^j$ on $U_h^j$ by the variational form in \eqref{eq:fd-rec}, with the integration band $U^{j-1}$ replaced by the pre-update band $U^{j,-}$. The semidiscrete scheme uses the new band $U[\phi_h(t_{m+1}^+)]$ here; lagging the band keeps this reconstruction linear. The accepted geometry is required to satisfy
$$
 \begin{gathered}
 \bigl\{x\in D:\dist(x,U^j)\le\tfrac12c_{\rm buf}h\bigr\}
 \Subset U_h^j
 \subset\bigl\{x\in D:\dist(x,U^j)\le2c_{\rm buf}h\bigr\},\\
 U_h^j\subset\Omega_h^j.
 \end{gathered}
$$

These conditions are checked after reconstruction on the new active mesh, together with nondegeneracy of the relevant level surfaces, consistency of the inherited signs, and $\inf_{x\in U_h^j}|\eta_h^j(x)|>0$. If a check fails, discard the trial transfer and mesh change and repeat from the preceding accepted state with a smaller $\tau_{j-1}$. Only accepted steps advance the time level. If no admissible step can be found, the computation is stopped; reducing the time step alone is not asserted to control a reconstruction jump for arbitrary $h$. A fully discrete analysis would also have to prove termination of these checks and a strict reset margin under its smallness assumptions.

Finally, given $\phi_h^j$ and $\Omega^j$, we compute
$u_h^j\in V_h|_{\Omega_h^j}$ from
\begin{equation}
 (u_h^j,v_h)_{\Omega^j}
 +(\nabla u_h^j,\nabla v_h)_{\Omega^j}
 +g_1(u_h^j,v_h)_{U_h^j}
 =(f^j,v_h)_{\Omega^j}
 \qquad
 \forall v_h\in V_h|_{\Omega_h^j},
 \label{eq:fd-u}
\end{equation}
where $f^j:=f(\cdot,t_j^{\rm fd})$. This problem is linear because $\Omega^j$ is known. The resulting fields $(p_h^j,\eta_h^j,\phi_h^j,u_h^j)$ provide the data for the next time step, and the active meshes are kept until \eqref{eq:fd-update-test} fails.

\section{The chain condition}
\label{appendix:path}
For $t\in I_m\cap[0,t^\ast)$ and sufficiently small $h$, the definitions in \eqref{eq:buffer-elements}--\eqref{eq:numerical-bulk-support}
ensure that every element $K$ in $\Uha$ contains a point in the $c_{\rm buf}h$-neighborhood of $U[\phi_h(t_m^-)]$. Joining this point to a closest point in $\overline{U[\phi_h(t_m^-)]}$, and preceding this segment by a segment inside $K$ if necessary, gives a path from the interior of $K$ to the reference band $U[\phi_h(t_m^-)]$ of length at most $Ch$ contained in $\Uha$, as in \cref{section:patch}. 

The endpoint of this path lies in $\overline{U[\phi_h(t_m^-)]}\subset\Uha$, so the upper-distance condition \eqref{eq:band-buffer-b} provides a point in the current band $\overline{U[\phi_h(t)]}$ at distance at most $2c_{\rm buf}h$ from it. This current-band point belongs to a selected mesh element. Applying the preceding construction to that element joins the current-band point to
$\overline{U[\phi_h(t_m^-)]}$ by another path of length at most $Ch$, contained in that element and the closed $c_{\rm buf}h$-neighborhood of $U[\phi_h(t_m^-)]$. The endpoints of these two paths in the reference band are therefore at distance at most $Ch$ from one another.

The uniform Lipschitz charts and fixed positive width of the reference band allow these two endpoints to be joined within
$\overline{U[\phi_h(t_m^-)]}$ by a polygonal path of length at most $Ch$. Indeed, for sufficiently small $h$, nearby boundary points lie in a common Lipschitz chart, where such a path is obtained by moving into the band, joining parallel to the base of the chart, and returning to the other endpoint. Away from the boundary, a straight segment suffices. Concatenating the paths, with the second one traversed in reverse, gives a path from the interior of $K$ to the current band, contained in $\Uha$ and of length at most $Ch$. The uniform interior-cone property of $U[\phi_h(t)]$ then allows this path to be continued, over an additional distance at most $Ch$, to a whole element contained in $U[\phi_h(t)]$ and at distance at least $h$ from its boundary.

These paths can be perturbed within $\Uha$ to cross relative interiors of mesh faces rather than vertices or lower-dimensional edges. The connecting portions lie inside selected elements, the closed $c_{\rm buf}h$-neighborhood of $U[\phi_h(t_m^-)]$, or the current band. Every element incident to a point of this neighborhood is selected by \eqref{eq:buffer-elements}, and the current band lies inside $\Uha$ with the clearance $\tfrac12c_{\rm buf}h$ supplied by \eqref{eq:band-buffer-a}. Thus the perturbations remain in $\Uha$, and all crossed faces belong to $\mathcal F_h^g(\Uha)$. Each resulting chain lies in a $Ch$-neighborhood of its starting element. Quasiuniformity and shape regularity therefore give a uniformly bounded number of face crossings, after removing repetitions, and uniformly bounded chain overlap. Hence the pair $(\mathcal B_{h,m},\Uh)$ satisfies the chain condition of \cref{def:chain}.

Since the boundaries of $U[\phi]$ and $U[\phi_h]$ are corresponding normal graphs with displacement at most $Ch^{1.1}$ and the same side classification, the terminal elements chosen at distance at least $h$ inside $U[\phi_h]$ also lie wholly inside $U[\phi]$ for sufficiently small $h$. The same chains therefore establish the chain condition for the pair $(\mathcal B_{h,m},U[\phi])$.

Moreover, the mesh construction and the interior/exterior classification outside $\Uha$, inherited from $\Omega[\phi_h(t_m^-)]$ and kept fixed on $I_m$, give
$$
 \Omega_h[\phi_h]\setminus\Omega[\phi_h]
 =
 U_h[\phi_h]\setminus\Omega[\phi_h].
$$
Indeed, every element of $\Omega_h[\phi_h]$ outside $U_h[\phi_h]$ lies wholly inside $\Omega[\phi_h]$. Thus only elements of
$U_h[\phi_h]$ require extension for the bulk estimate. Apply the preceding construction to each such element. If its terminal element is not already contained in $\Omega[\phi_h]$, the exterior width $\kappa h$ of $U[\phi_h]$ and the uniform normal-graph geometry allow the path to be continued inside the band, over a further distance at most $Ch$, to a whole element contained in $\Omega[\phi_h]\cap U[\phi_h]$. All crossed faces remain in $\mathcal F_h^g(\Uha)$, and the same local element-count argument gives uniformly bounded chain length and overlap. Hence the pair $(\T_{h,m}^\Omega,\Oh)$ also satisfies the chain condition of \cref{def:chain}. The same argument applies at $t_{m+1}^-$ on the old meshes, by the non-strict version of \eqref{eq:band-buffer}, and at $t_{m+1}^+$ on the new meshes, where $U[\phi_h(t_{m+1}^-)]$ is the reference band.

%The same argument applies to the old-support left limits at update times, using the limiting inclusions in \eqref{eq:band-buffer}.

\section{Completion of the proof}
\label{sec:closure}
All constants in the spatial estimates are uniform in the update index and the cut position. Constants used only for local finite-dimensional existence may depend on the fixed value of $h$.

\subsection{Closing the continuation argument}
\label{sec:closing}

In \cref{sec:ODE} we show the following results: (i) the semidiscrete problem can be regarded as an ODE system for $(p_h,\eta_h)$ on every time interval $I_m$ where the active mesh is fixed; (ii) at each time, for transport data with small errors, equation \eqref{eq:scheme-rec} has a unique solution on the local branch, and this solution satisfies \eqref{eq:ephi-final} with $\Uh$ replaced by $\Uha$ on the right-hand side. 

Since $e_p=0$ and $e_\eta=0$ at time $t=0$, the local existence, uniqueness and continuity of ODE systems mentioned in (i) imply that the semidiscrete problem has a solution satisfying $\nL{e_p}{\Uha}+\nL{e_\eta}{\Uha}\leq h^{k}$ up to some positive time, and \eqref{eq:eu-final}--\eqref{eq:ephi-final} imply that $\nH{e_\phi}{\Uha}+\nH{e_u}{\Oha}\le Ch^{k}$ up to this positive time. Thus, for sufficiently small $h$, the continuation conditions \eqref{eq:continuation-conditions} hold up to some positive time. This implies that $t^\ast>0$. Moreover, by \eqref{time-length-min} only finitely many updates occur before $t^\ast$, so the active meshes are fixed on some interval $(t_m,t^\ast)$. The solution exists up to $t^\ast$ itself because the nodal vector associated with the finite element solution stays in the compact set where \eqref{eq:continuation-conditions} holds, on which $|\eta_h|$ is bounded away from zero by \eqref{eq:smallness}. There the right-hand side of the ODE is bounded and Lipschitz, so the solution is uniformly continuous and extends to $t^\ast$. The conditions \eqref{eq:continuation-conditions} then hold at $t^\ast$ because they are closed. 

If $t^\ast<T$ is not an update time, the errors at $t^\ast$ are bounded by $C_0h^k$ by \cref{prop:coupled} and continuity, and therefore bounded by $\frac12 h^{k-0.1}$ for sufficiently small $h$. The local existence and uniqueness of solutions to the ODE system imply that the solution exists up to a time bigger than $t^\ast$, and the continuity in time on the fixed active mesh ensures that \eqref{eq:continuation-conditions} is satisfied on $[t^\ast,t^\ast+\epsilon^\ast]$ for some $\epsilon^\ast>0$; this contradicts the definition of $t^\ast$. 

If $t^\ast<T$ is an update time, the errors at $t^{\ast,-}$ are bounded by $C_0h^k$ according to \cref{prop:coupled}, and the errors after being transferred to the new active mesh are still bounded by $Ch^k$ according to \cref{lem:transfer}. Hence, the errors at $t^{\ast,+}$ are bounded by $\frac12 h^{k-0.1}$ for sufficiently small $h$. As in the previous case, the solution on the new meshes satisfies \eqref{eq:continuation-conditions} on $[t^\ast,t^\ast+\epsilon^\ast]$ for some $\epsilon^\ast>0$, which again contradicts the definition of $t^\ast$. 

Hence $t^\ast=T$ for sufficiently small $h$. Passing to the limit $t\to T$ gives \eqref{eq:main-theorem-error} on $[0,T]$, and \eqref{time-length-min} with $t^\ast=T$ gives \eqref{eq:event-time interval-spacing}. This proves \cref{thm:main}.

\subsection{Semidiscrete problem as an ODE system}
\label{sec:ODE}

We prove the two properties stated in \cref{sec:closing}. To justify the equivalence of the semidiscrete problem to an ODE system for $(p_h,\eta_h)$, we show that, at each fixed time, $\phi_h$ is locally uniquely determined by $p_h$ and $\eta_h$ and depends locally Lipschitz continuously on these data. The bulk equation \eqref{eq:scheme-u} then determines $u_h$, and the transport equations \eqref{eq:scheme-p}--\eqref{eq:scheme-eta} determine $(\partial_tp_h,\partial_t\eta_h)$ through their invertible stabilized mass matrix. This determines $(\partial_tp_h,\partial_t\eta_h)$ as a function of $(p_h,\eta_h)$. Hence the semidiscrete problem can be regarded as an ODE system for $(p_h,\eta_h)$.

Since $\phi_h=\phi+\rho_\phi+e_\phi$, with $\phi$ and $\rho_\phi$ determined by the exact solution, solving \eqref{eq:scheme-rec} for $\phi_h$ is equivalent to solving \eqref{eq:ephi-current-band} for $e_\phi$. We use the following fixed-point argument.

Fix $t\in[t_m,t_{m+1}]$, with the values $t_m^+$ and $t_{m+1}^-$ at the endpoints; on this interval the active meshes, built from $\phi_h(t_m^-)$, do not change. For given transport data satisfying
$$
 \|e_p\|_{L^2(U_h[\phi_h(t_m^-)])}+\|e_\eta\|_{L^2(U_h[\phi_h(t_m^-)])}
 \le 2h^{k-0.1},
$$
set
$$
 \mathcal B_h
 :=
 \bigl\{\chi\in V_h|_{U_h[\phi_h(t_m^-)]}:
        \|\chi\|_{L^\infty(U_h[\phi_h(t_m^-)])}\le h^{1.1}\bigr\}.
$$
For $\chi\in\mathcal B_h$, the finite element inverse estimate gives $\|\nabla\chi\|_{L^\infty(U_h[\phi_h(t_m^-)])}\le Ch^{0.1}$. This shows that, for sufficiently small $h$, the band $U[\phi_h^\ast+\chi]$ is nondegenerate and remains inside the active mesh $U_h[\phi_h(t_m^-)]$. We use the geometric and ghost-extension estimates also for these nearby trial bands.

For $\chi\in\mathcal B_h$, define $\mathcal M(\chi)\in V_h|_{U_h[\phi_h(t_m^-)]}$ to be the finite element solution of the following weak formulation: 
\begin{equation}
\begin{aligned}
 &(\mathcal M(\chi),v_h)_{U[\phi_h^\ast+\chi]}
 +(\nabla\mathcal M(\chi),\nabla v_h)_{U[\phi_h^\ast+\chi]}
 +g_1(\mathcal M(\chi),v_h)_{U_h[\phi_h(t_m^-)]}\\
 &\quad=(e_p,v_h)_{U[\phi_h^\ast+\chi]}
 +(e_\eta,\nabla v_h)_{U[\phi_h^\ast+\chi]}\\
 &\qquad
 +(\rho_\eta-\nabla\rho_\phi,\nabla v_h)_{U[\phi_h^\ast+\chi]}
 -g_1(\rho_\phi,v_h)_{U_h[\phi_h(t_m^-)]}
 \qquad\forall v_h\in V_h|_{U_h[\phi_h(t_m^-)]},
\end{aligned}
 \label{eq:ephi-current-band-chi}
\end{equation}
where $\phi_h^\ast=I_h\phi(t)$. The integration band is fixed for a given $\chi$, and the stabilized bilinear form is coercive
by the ghost-extension estimate. Thus, \eqref{eq:ephi-current-band-chi} has a unique solution. Taking $v_h=\mathcal M(\chi)$ and using this coercivity, \eqref{eq:ephi-interp-residual} and Young's inequality yields
$$
 \|\mathcal M(\chi)\|_{H^1(U_h[\phi_h(t_m^-)])}^2
 \le
 C\|e_p\|_{L^2(U_h[\phi_h(t_m^-)])}^2
 +C\|e_\eta\|_{L^2(U_h[\phi_h(t_m^-)])}^2
 +Ch^{2k}
 \le Ch^{2k-0.2}.
$$
The discrete Sobolev inequality in \cref{lem:discrete-sobolev} consequently gives
$$
 \|\mathcal M(\chi)\|_{L^\infty(U_h[\phi_h(t_m^-)])}
 \le Ch^{k-0.1-(n-1)/4}
 \le \tfrac12h^{1.1}
$$
for sufficiently small $h$. Hence $\mathcal M:\mathcal B_h\to\mathcal B_h$ is well defined.

It remains to show that $\mathcal M$ is a contraction. Let $U_i:=U[\phi_h^\ast+\chi_i]$, $i=1,2$. Subtracting the two equations \eqref{eq:ephi-current-band-chi}, we obtain
\begin{equation}
\begin{aligned}
 &(\mathcal M(\chi_1)-\mathcal M(\chi_2),v_h)_{U_1}\\
 &\quad
 +(\nabla(\mathcal M(\chi_1)-\mathcal M(\chi_2)),
      \nabla v_h)_{U_1}
 +g_1(\mathcal M(\chi_1)-\mathcal M(\chi_2),v_h)_{U_h[\phi_h(t_m^-)]}\\
 &=(e_p-\mathcal M(\chi_2),v_h)_{U_1}
   -(e_p-\mathcal M(\chi_2),v_h)_{U_2}\\
 &\quad
 +(e_\eta-\nabla\mathcal M(\chi_2),\nabla v_h)_{U_1}
 -(e_\eta-\nabla\mathcal M(\chi_2),\nabla v_h)_{U_2}\\
 &\quad
 +(\rho_\eta-\nabla\rho_\phi,\nabla v_h)_{U_1}
 -(\rho_\eta-\nabla\rho_\phi,\nabla v_h)_{U_2}
 \qquad\forall v_h\in V_h|_{U_h[\phi_h(t_m^-)]}.
\end{aligned}
 \label{eq:ephi-current-band-chi12}
\end{equation}

Each difference on the right-hand side is estimated using \eqref{eq:FTC-band}. For example, let $\chi^\theta=(1-\theta)\chi_2+\theta\chi_1$. The second term on the right-hand side of \eqref{eq:ephi-current-band-chi12} can be bounded by
\begin{align*}
 &C\int_0^1
 \|e_\eta-\nabla\mathcal M(\chi_2)\|
      _{L^2(\Gamma_\pm[\phi_h^\ast+\chi^\theta])}
 \|\nabla v_h\|_{L^2(\Gamma_\pm[\phi_h^\ast+\chi^\theta])}
 \|\chi_1-\chi_2\|_{L^\infty(U_h[\phi_h(t_m^-)])}
 \,\dd\theta\\
 &\quad\le
 Ch^{-1}
\|e_\eta-\nabla\mathcal M(\chi_2)\|_{L^2(U_h[\phi_h(t_m^-)])}
 \|v_h\|_{H^1(U_h[\phi_h(t_m^-)])}
 \|\chi_1-\chi_2\|_{L^\infty(U_h[\phi_h(t_m^-)])}\\
 &\quad\le
 Ch^{k-1.1}\|v_h\|_{H^1(U_h[\phi_h(t_m^-)])}
 \|\chi_1-\chi_2\|_{L^\infty(U_h[\phi_h(t_m^-)])},
\end{align*}
where the second inequality uses the elementwise inverse trace estimate. The first term on the right-hand side of \eqref{eq:ephi-current-band-chi12} has the same bound. For the third term on the right-hand side of \eqref{eq:ephi-current-band-chi12}, $\|\rho_\eta-\nabla\rho_\phi\|_{L^\infty}\le Ch^k$ and one inverse trace estimate give the bound
$$
 Ch^{k-1/2}\|v_h\|_{H^1(U_h[\phi_h(t_m^-)])}
 \|\chi_1-\chi_2\|_{L^\infty(U_h[\phi_h(t_m^-)])} .
$$
Therefore, taking $v_h=\mathcal M(\chi_1)-\mathcal M(\chi_2)$ in \eqref{eq:ephi-current-band-chi12} and using coercivity, we obtain
$$
 \|\mathcal M(\chi_1)-\mathcal M(\chi_2)\|_{H^1(U_h[\phi_h(t_m^-)])}
 \le Ch^{k-1.1}
       \|\chi_1-\chi_2\|_{L^\infty(U_h[\phi_h(t_m^-)])}.
$$
The discrete Sobolev inequality in \cref{lem:discrete-sobolev} therefore yields
$$
 \|\mathcal M(\chi_1)-\mathcal M(\chi_2)\|_{L^\infty(U_h[\phi_h(t_m^-)])}
 \le Ch^{k-1.1-(n-1)/4}
       \|\chi_1-\chi_2\|_{L^\infty(U_h[\phi_h(t_m^-)])}.
$$
Since $k-1.1-(n-1)/4\ge1.4$ for $k\ge3$ and $n\le3$, the contraction factor is at most $1/2$ for sufficiently small $h$. Banach's fixed-point theorem gives a unique fixed point $e_\phi\in\mathcal B_h$, and hence a locally unique solution $\phi_h=\phi_h^\ast+e_\phi$ of \eqref{eq:scheme-rec}. 

Once $\phi_h$ has been determined, $u_h$ can be determined by \eqref{eq:scheme-u}. The local Lipschitz dependence of $u_h$ and $\phi_h$ on $(p_h,\eta_h)$ is proved in the same way.


\begin{thebibliography}{99}

\bibitem{BadiaCodina2006}
S.~Badia and R.~Codina,
\newblock Analysis of a stabilized finite element approximation of the transient convection--diffusion equation using an ALE framework, 
\newblock \emph{SIAM J. Numer. Anal.}, 44 (2006), pp.~2159--2197. %\url{https://doi.org/10.1137/050643532}.

\bibitem{BarrettGarckeNurnberg2010}
J.~W. Barrett, H.~Garcke, and R.~N\"urnberg,
\newblock On stable parametric finite element methods for the Stefan problem and the Mullins--Sekerka problem with applications to dendritic growth,
\newblock \emph{J. Comput. Phys.}, 229 (2010), pp.~6270--6299. 
%\url{https://doi.org/10.1016/j.jcp.2010.04.039}.

\bibitem{BrennerScott2008}
S.~C. Brenner and L.~R. Scott, \emph{The Mathematical Theory of Finite
Element Methods}, 3rd ed., Springer, New York, 2008.

\bibitem{Burman2010}
E.~Burman,
\newblock Ghost penalty,
\newblock \emph{Comptes Rendus. Mathématique}, 348 (2010), pp.~1217--1220. 
%\url{https://doi.org/10.1016/j.crma.2010.10.006}.

\bibitem{BCDE2021}
E.~Burman, M.~Cicuttin, G.~Delay, and A.~Ern,
\newblock An unfitted hybrid high-order method with cell agglomeration for elliptic interface problems,
\newblock \emph{SIAM J. Sci. Comput.}, 43 (2021), pp.~A859--A882.

\bibitem{BurmanEtAl2015}
E.~Burman, S.~Claus, P.~Hansbo, M.~G.~Larson, and A.~Massing,
\newblock CutFEM: discretizing geometry and partial differential equations,
\newblock \emph{Internat. J. Numer. Methods Engrg.}, 104 (2015), pp.~472--501.

\bibitem{BurmanFernandez2014}
E.~Burman and M.~A. Fern\'andez,
\newblock An unfitted Nitsche method for incompressible fluid-structure interaction using overlapping meshes,
\newblock \emph{Comput. Methods Appl. Mech. Engrg.}, 279 (2014), pp.~497--514.

\bibitem{BurmanFreiMassing2022}
E.~Burman, S.~Frei, and A.~Massing,
\newblock Eulerian time-stepping schemes for the non-stationary Stokes
equations on time-dependent domains,
\newblock \emph{Numer. Math.}, 150 (2022), pp.~423--478.

\bibitem{BurmanHansboLarson2022}
E.~Burman, P.~Hansbo, and M.~G.~Larson,
\newblock CutFEM based on extended finite element spaces,
\newblock \emph{Numer. Math.}, 152 (2022), pp.~331--369.

\bibitem{Ciarlet1978}
P.~G. Ciarlet, \emph{The Finite Element Method for Elliptic Problems},
North-Holland, Amsterdam, 1978.

\bibitem{ClausKerfriden2019} S.~Claus and P.~Kerfriden, \newblock A CutFEM method for two-phase flow problems, \newblock \emph{Comput. Methods Appl. Mech. Engrg.}, 348 (2019), pp.~185--206.

\bibitem{DeckelnickEtAl2010}
K.~Deckelnick, G.~Dziuk, C.~M. Elliott, and C.-J. Heine,
\newblock An $h$-narrow band finite-element method for elliptic equations on
implicit surfaces,
\newblock \emph{IMA J. Numer. Anal.}, 30 (2010), pp.~351--376.

\bibitem{DeckelnickElliottRanner2014}
K.~Deckelnick, C.~M. Elliott, and T.~Ranner,
\newblock Unfitted finite element methods using bulk meshes for surface
partial differential equations,
\newblock \emph{SIAM J. Numer. Anal.}, 52 (2014), pp.~2137--2162.

\bibitem{DelfourZolesio2011}
M.~C. Delfour and J.-P. Zol\'esio, \emph{Shapes and Geometries: Metrics,
Analysis, Differential Calculus, and Optimization}, 2nd ed., SIAM,
Philadelphia, 2011.

\bibitem{DziukElliott2007}
G.~Dziuk and C.~M. Elliott,
\newblock Finite elements on evolving surfaces,
\newblock \emph{IMA J. Numer. Anal.}, 27 (2007), pp.~262--292.
%\url{https://doi.org/10.1093/imanum/drl023}.

\bibitem{DziukElliott2013}
G.~Dziuk and C.~M. Elliott,
\newblock Finite element methods for surface PDEs,
\newblock \emph{Acta Numer.}, 22 (2013), pp.~289--396.

\bibitem{EdelmannKovacsLubich2025}
D.~Edelmann, B.~Kov\'acs, and C.~Lubich,
\newblock Numerical analysis of an evolving bulk-surface model of tumour growth,
\newblock \emph{IMA J. Numer. Anal.}, 45 (2025), pp.~2581--2627.
%\url{https://doi.org/10.1093/imanum/drae077}.

\bibitem{ElliottRanner2021}
C.~M. Elliott and T.~Ranner,
\newblock A unified theory for continuous-in-time evolving finite element
space approximations to partial differential equations in evolving domains,
\newblock \emph{IMA J. Numer. Anal.}, 41 (2021), pp.~1696--1845.
%\url{https://doi.org/10.1093/imanum/draa062}.

\bibitem{EylesKingStyles2019}
J. Eyles, J. R. King, and V. Styles,
\newblock A tractable mathematical model for tissue growth,
\newblock \emph{Interfaces Free Bound.}, 21 (2019), pp.~463--493. 
%\url{https://doi.org/10.4171/IFB/428}.

\bibitem{FernandezLiOlshanskii2026}
M.~A. Fern\'andez, B.~Li, and M.~A. Olshanskii,
\newblock Convergence of a CutFEM for fluid--structure interaction with a
deforming interface,
\newblock \url{https://doi.org/10.48550/arXiv.2608.08140}, 2026.

\bibitem{FrachonZahedi2019} 
T.~Frachon and S.~Zahedi, \newblock A cut finite element method for incompressible two-phase Navier--Stokes flows, \newblock \emph{J. Comput. Phys.}, 384 (2019), pp.~77--98.

\bibitem{FriedmanReitich1999}
A.~Friedman and F.~Reitich,
\newblock Analysis of a mathematical model for the growth of tumors,
\newblock \emph{J. Math. Biol.}, 38 (1999), pp.~262--284. 
%\url{https://doi.org/10.1007/s002850050149}.

\bibitem{FriesBelytschko2010}
T.-P. Fries and T.~Belytschko,
\newblock The extended/generalized finite element method: an overview of the
method and its applications,
\newblock \emph{Internat. J. Numer. Methods Engrg.}, 84 (2010), pp.~253--304.

\bibitem{GawlikLew2015}
E.~S. Gawlik and A.~J. Lew,
\newblock Unified analysis of finite element methods for problems with
moving boundaries,
\newblock \emph{SIAM J. Numer. Anal.}, 53 (2015), pp.~2822--2846. 
%\url{https://doi.org/10.1137/140990437}.

\bibitem{GrossReusken2007}
S.~Gro{\ss} and A.~Reusken,
\newblock An extended pressure finite element space for two-phase
incompressible flows with surface tension,
\newblock \emph{J. Comput. Phys.}, 224 (2007), pp.~40--58.

%\bibitem{GrossReusken2011}
%S.~Gro{\ss} and A.~Reusken,
%\newblock \emph{Numerical Methods for Two-phase Incompressible Flows},
%Springer Ser. Comput. Math. 40, Springer, Berlin, 2011.

\bibitem{HansboHansbo2002}
A.~Hansbo and P.~Hansbo,
\newblock An unfitted finite element method, based on Nitsche's method, for
elliptic interface problems,
\newblock \emph{Comput. Methods Appl. Mech. Engrg.}, 191 (2002), pp.~5537--5552.

\bibitem{HansboLarsonZahedi2016}
P.~Hansbo, M.~G. Larson, and S.~Zahedi,
\newblock A cut finite element method for coupled bulk-surface problems on
time-dependent domains,
\newblock \emph{Comput. Methods Appl. Mech. Engrg.}, 307 (2016), pp.~96--116.

%\bibitem{Hu-Li-2022} J.~Hu and B.~Li, \newblock Evolving finite element methods with an artificial tangential velocity for mean curvature flow and Willmore flow, \newblock \emph{Numer. Math.}, 152 (2022), no.~1, pp.~127--181.

\bibitem{JohanssonLarson2013} 
A.~Johansson and M.~G.~Larson, \newblock A high order discontinuous Galerkin Nitsche method for elliptic problems with fictitious boundary, \newblock \emph{Numer. Math.}, 123 (2013), no.~4, pp.~607--628.

\bibitem{KLLP2017} 
B.~Kovács, B.~Li, C.~Lubich and C.~A.~Power Guerra, 
\newblock Convergence of finite elements on an evolving surface driven by diffusion on the surface, \newblock \emph{Numer. Math.}, 137 (2017), no.~3, pp.~643--689.

\bibitem{KLL2019}
B.~Kovács, B.~Li and C.~Lubich,
\newblock A convergent evolving finite element algorithm for mean curvature flow of closed surfaces,
\newblock \emph{Numer. Math.}, 143 (2019), no.~4, pp.~797--853.

\bibitem{KLL2021}
B.~Kovács, B.~Li and C.~Lubich,
\newblock A convergent evolving finite element algorithm for Willmore flow of closed surfaces,
\newblock \emph{Numer. Math.}, 149 (2021), no.~3, pp.~595--643.

\bibitem{Lehrenfeld2016}
C.~Lehrenfeld,
\newblock High order unfitted finite element methods on level set domains
using isoparametric mappings,
\newblock \emph{Comput. Methods Appl. Mech. Engrg.}, 300 (2016), pp.~716--733.

\bibitem{LehrenfeldOlshanskii2019}
C.~Lehrenfeld and M.~A. Olshanskii,
\newblock An Eulerian finite element method for PDEs in time-dependent domains,
\newblock \emph{ESAIM Math. Model. Numer. Anal.}, 53 (2019), pp.~585--614.

\bibitem{LehrenfeldOlshanskiiXu2018}
C.~Lehrenfeld, M.~A. Olshanskii, and X.~Xu,
\newblock A stabilized trace finite element method for partial differential
equations on evolving surfaces,
\newblock \emph{SIAM J. Numer. Anal.}, 56 (2018), pp.~1643--1672.

\bibitem{LehrenfeldReusken2018}
C.~Lehrenfeld and A.~Reusken,
\newblock Analysis of a high-order unfitted finite element method for
elliptic interface problems,
\newblock \emph{IMA J. Numer. Anal.}, 38 (2018), pp.~1351--1387.

\bibitem{LouLehrenfeld2022}
Y.~Lou and C.~Lehrenfeld,
\newblock Isoparametric unfitted BDF--{Finite Element Method} for PDEs on
evolving domains,
\newblock \emph{SIAM J. Numer. Anal.}, 60 (2022), pp.~2069--2098. 
%\url{https://doi.org/10.1137/21M142126X}.

%\bibitem{MassingEtAl2014}
%A.~Massing, M.~G. Larson, A.~Logg, and M.~E. Rognes,
%\newblock A stabilized Nitsche fictitious domain method for the Stokes
%problem,
%\newblock \emph{J. Sci. Comput.}, 61 (2014), pp.~604--628.

\bibitem{OlshanskiiReusken2026}
M.~A. Olshanskii and A.~Reusken,
\newblock Analysis of a finite element method for PDEs in evolving domains
with topological changes,
\newblock \emph{Interfaces Free Bound.}, published online June~10, 2026,
\url{https://doi.org/10.4171/IFB/572}.

\bibitem{OlshanskiiReuskenGrande2009}
M.~A. Olshanskii, A.~Reusken, and J.~Grande,
\newblock A finite element method for elliptic equations on surfaces,
\newblock \emph{SIAM J. Numer. Anal.}, 47 (2009), pp.~3339--3358.

\bibitem{OlshanskiiReuskenXu2014}
M.~A. Olshanskii, A.~Reusken, and X.~Xu,
\newblock An Eulerian space-time finite element method for diffusion
problems on evolving surfaces,
\newblock \emph{SIAM J. Numer. Anal.}, 52 (2014), pp.~1354--1377.

\bibitem{Saye2015}
R.~I. Saye,
\newblock High-order quadrature methods for implicitly defined surfaces and volumes in hyperrectangles,
\newblock \emph{SIAM J. Sci. Comput.}, 37 (2015), pp.~A993--A1019.

\bibitem{TchindaNgueyongUrquizaMartin2024a} 
I.~Tchinda Ngueyong, J.~M.~Urquiza and D.~Martin, \newblock A CutFEM method for phase change problems with natural convection, \newblock \emph{Comput. Methods Appl. Mech. Engrg.}, 420 (2024), Paper No.~116713, 24~pp.

\bibitem{TchindaNgueyongUrquizaMartin2024b} 
I.~Tchinda Ngueyong, J.~M.~Urquiza and D.~Martin, \newblock A Nitsche-based cut finite element solver for two-phase Stefan problems, \newblock \emph{Internat. J. Numer. Methods Engrg.}, 125 (2024), no.~6, Paper No.~e7408, 25~pp.

\bibitem{vanNoordenEck2011}
T.~L. van Noorden and C.~Eck,
\newblock Phase field approximation of a kinetic moving-boundary problem
modelling dissolution and precipitation,
\newblock \emph{Interfaces Free Bound.}, 13 (2011), pp.~29--55. 
%\url{https://doi.org/10.4171/IFB/247}.

\bibitem{vonWahlRichter2023}
H.~von Wahl and T.~Richter,
\newblock Error analysis for a parabolic PDE model problem on a coupled
moving domain in a fully Eulerian framework,
\newblock \emph{SIAM J. Numer. Anal.}, 61 (2023), pp.~286--314.

\bibitem{vonWahlRichterLehrenfeld2022}
H.~von Wahl, T.~Richter, and C.~Lehrenfeld,
\newblock An unfitted Eulerian finite element method for the time-dependent Stokes problem on moving domains,
\newblock \emph{IMA J. Numer. Anal.}, 42 (2022), pp.~2505--2544.

\bibitem{Zahedi2018}
S.~Zahedi,
\newblock A space-time cut finite element method with quadrature in time,
\newblock in pp.~281--306 of \emph{Geometrically Unfitted Finite Element Methods and Applications}, edited by S.~P.~A. Bordas, E.~Burman, M.~G. Larson, and M.~A. Olshanskii. Lecture Notes in Computational Science and Engineering, vol 121. (2017), Springer, Cham. 

\end{thebibliography}
\end{document}